\documentclass{amsart}

\usepackage{amsmath}
\usepackage[latin1]{inputenc}
\usepackage{amsfonts}
\usepackage{amssymb}
\usepackage{graphicx,epsfig}
\usepackage{amsthm}
\usepackage{amscd}
\usepackage[all]{xy}
\usepackage{verbatim}
\usepackage{pst-all}

\usepackage{amsmath}

\usepackage{amsfonts}

\usepackage{amssymb}

\usepackage{amsthm}

\usepackage{graphics}

\newcommand{\mk}{\medskip}

\newcommand{\ZZ}{\mathbb{Z}}

\newcommand{\CC}{\mathbb{C}}

\newcommand{\NN}{\mathbb{N}}

\newcommand{\QQ}{\mathbb{Q}}

\newcommand{\Glie}{\mathfrak{g}}

\newcommand{\Yim}{\mathcal{Y}}

\newcommand{\Hlie}{\mathfrak{h}}

\newcommand{\demo}{\noindent {\it \small Proof:}\quad}

\renewcommand{\NN}{\ensuremath{\mathbb{N}}}

\renewcommand{\CC}{\ensuremath{\mathbb{C}}}

\renewcommand{\QQ}{\ensuremath{\mathbb{Q}}}

\newcommand{\U}{\mathcal{U}}

\newtheorem{thm}{Theorem}[section]

\newtheorem{defi}[thm]{Definition}

\newtheorem{cor}[thm]{Corollary}

\newtheorem{prop}[thm]{Proposition}

\newtheorem{lem}[thm]{Lemma}

\newtheorem{rem}[thm]{Remark}

\author{David Hernandez}
\address{\'Ecole Normale Sup\'erieure Paris and
Univ. di Roma La Sapienza
\newline Address : ENS-DMA, 45, Rue d'Ulm F-75230 PARIS,
  Cedex 05  FRANCE}
\email{David.Hernandez@ens.fr\\ URL: http://www.dma.ens.fr/$\sim$dhernand}

\usepackage[all]{xy}
\title{Drinfeld coproduct, quantum fusion tensor category and applications}

\begin{document}

\begin{abstract} The class of quantum affinizations (or quantum loop algebras, see
\cite{Dri2, Cha2, gkv, var, mi, Naams, jin, her04}) includes quantum affine algebras
and quantum toroidal algebras. In general they have no Hopf
algebra structure, but have a ``coproduct'' (the Drinfeld coproduct) which
does not produce tensor products of modules in the usual way because
it is defined in a completion. In this paper we propose a new
process to produce quantum fusion modules from it : for all quantum
affinizations, we construct by deformation and renormalization a new
(non semi-simple) tensor category Mod. For quantum affine algebras
this process is new and different from the usual tensor product. For general
quantum affinizations, for example for toroidal algebras, so far, no process
to produce fusion modules was known. We derive several applications from
it :  we construct the fusion of (finitely many) arbitrary $l$-highest 
weight modules, and
prove that it is always cyclic. We establish exact sequences
involving fusion of Kirillov-Reshetikhin modules related to new
$T$-systems extending results of \cite{Nab, Nad, her06}. Eventually for a large
class of quantum affinizations we prove that the subcategory of
finite length modules of Mod is stable under the new monoidal
bifunctor.

\vskip 4.5mm


\end{abstract}

\maketitle

\tableofcontents

\section{Introduction} 

In this paper $q\in\mathbb{C}^*$ is not a root of unity.

\noindent Drinfeld \cite{Dri1} and Jimbo \cite{jim} associated, independently, to any symmetrizable Kac-Moody algebra $\mathfrak{g}$ and $q\in\mathbb{C}^*$ a Hopf algebra $\mathcal{U}_q(\mathfrak{g})$ called quantum Kac-Moody algebra. The quantum algebras of finite type $\mathcal{U}_q(\mathfrak{g})$ ($\mathfrak{g}$ of finite type) and their representations have been intensively studied (see for example \cite{Cha2, lu, ro} and references therein). The quantum affine algebras $\mathcal{U}_q(\hat{\mathfrak{g}})$ ($\hat{\mathfrak{g}}$ affine algebra) are also of particular interest : they have two realizations, the usual Drinfeld-Jimbo realization and a new realization (see \cite{Dri2, bec}) as a quantum affinization of a quantum algebra of finite type $\mathcal{U}_q(\mathfrak{g})$. Quantum affine algebras and their representations have also been intensively studied (see among others \cite{Aka, Cha0, Cha, Cha2, em, Fre, Fre2, her06, kas, Naams, Nab, var3} and references therein).

The quantum affinization process (that Drinfeld \cite{Dri2} described for constructing the second realization of a quantum affine algebra) can be extended to all symmetrizable quantum Kac-Moody algebras $\mathcal{U}_q(\mathfrak{g})$ (see \cite{jin, Naams}). One obtains a new class of algebras called quantum affinizations : the quantum affinization of $\mathcal{U}_q(\mathfrak{g})$ is denoted by $\mathcal{U}_q(\hat{\mathfrak{g}})$ and contains $\U_q(\Glie)$ as a subalgebra. The quantum affine algebras are the simplest examples and have the singular property of being also quantum Kac-Moody algebras. The quantum toroidal algebras (affinizations of quantum affine algebras) are also of
particular importance, see for example \cite{gkv, mi, mi2, Naams,
Nac, Sa, Sc, stu, tu, var}. In analogy to the Frobenius-Schur-Weyl duality between quantum groups of finite type and Hecke algebras \cite{jim2}, and between quantum affine algebras and affine Hecke algebras \cite{cp4}, quantum toroidal algebras have a close relation \cite{var1} to double affine
Hecke algebras and their degenerations (Cherednik's algebras \cite{che}) which have been recently
intensively studied (for example see \cite{be, che1, ggor, gs, va1, var2} and references therein). In general a quantum affinization (for example a quantum toroidal algebra) is not isomorphic to a quantum
Kac-Moody algebra and has no Hopf algebra structure. For convenience of the reader, we "describe" in the following diagram the relations between these different algebras (in particular with the two main classes of algebras considered in this paper, quantum Kac-Moody algebras on the left and quantum affinizations on the right).

\begin{center}
\epsfig{file=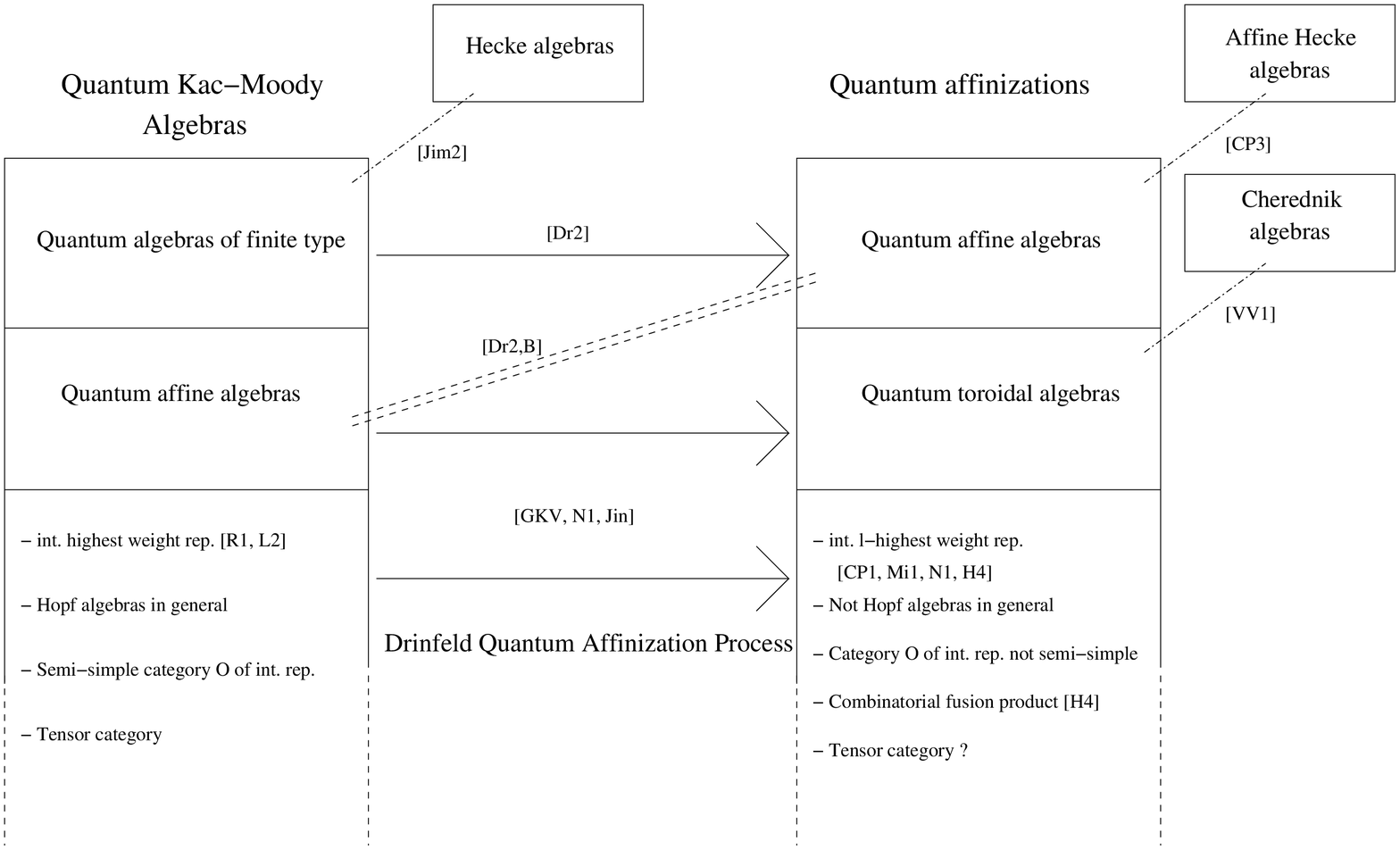,width=1\linewidth}
\end{center}

\noindent We denote by $\Hlie$ (resp. $\U_q(\Hlie)$) the Cartan subalgebra of $\Glie$ (resp. $\U_q(\Glie)$). $\U_q(\hat{\Hlie})\subset \U_q(\hat{\Glie})$ is the quantum affine analog of the Cartan subalgebra.

\noindent In analogy to the generalization of the representation
theory of quantum algebras of finite type to general quantum
Kac-Moody algebras (see \cite{lu}), a natural question is to extend
the representation theory of quantum affine algebras to general
quantum affinizations. By extending the results of \cite{Cha2}
(quantum affine algebras), \cite{mi} (quantum toroidal algebras of
type $A$) and \cite{Naams} (general simply-laced quantum
affinizations), we established in \cite{her04} a triangular
decomposition of $\U_q(\hat{\Glie})$ which allowed us to begin the
study of the representation theory of general quantum affinizations
: one can define the notion of $l$-highest weight representations of
$\U_q(\hat{\Glie})$ (analog to the usual notion of highest weight
module corresponding to $\U_q(\hat{\Hlie})$). One says that a $\U_q(\hat{\Glie})$-module is integrable 
(resp. in the category $\mathcal{O}$) if it is integrable (resp. in the category $\mathcal{O}$) as a $\U_q(\Glie)$-module.
Let $\text{Mod}(\U_q(\hat{\Glie}))$ be the category of integrable $\U_q(\hat{\Glie})$-module which are in the category $\mathcal{O}$.
It appears that $\U_q(\hat{\Glie})$ has numerous $l$-highest weight
representations in $\text{Mod}(\U_q(\hat{\Glie}))$ and that one can classify the simple ones. The category $\text{Mod}(\U_q(\hat{\Glie}))$ is not semi-simple. 

No Hopf algebra structure is known for $\U_q(\hat{\Glie})$ except in the
case of quantum affine algebras. However in \cite{her04} we proved
via Frenkel-Reshetikhin $q$-characters the existence of a
combinatorial fusion product in the Grothendieck group of
$\text{Mod}(\U_q(\hat{\Glie}))$, that it to say a product in the
semi-simplified category. As $\text{Mod}(\U_q(\hat{\Glie}))$ is not
semi-simple, this fusion product can not a priori be directly
translated in terms of modules.

\noindent Each quantum affinization has a "coproduct" (the Drinfeld
coproduct) which is defined in a completion of
$\U_q(\hat{\Glie})\otimes\U_q(\hat{\Glie})$. Although it can not be
directly used to construct tensor products of integrable
representations (it involves infinite sums), we propose in the
present paper a new process to produce quantum fusion modules from
it : for all quantum affinizations we construct by deformation of the Drinfeld coproduct and
renormalization a new (non semi-simple) tensor category Mod. In
particular it gives a representation theoretical interpretation of
the combinatorial fusion product and we get several applications from it. For
quantum affine algebras this process is new and different from the
usual tensor product (we prove that it has different properties),
and for general quantum affinizations (as general quantum toroidal
algebras) no process to produce fusion of modules was known so far.

\noindent For this construction, the first technical point that we
solve is a rationality problem for the dependence in the deformation
parameter $u$ of the Drinfeld coproduct. Let $\U_q'(\hat{\Glie})=\U_q(\hat{\Glie})\otimes\CC(u)$. We consider a "good" category $\text{Mod}(\U_q'(\hat{\Glie}))$ of $\U_q'(\hat{\Glie})$-modules $V$ with an integrable $\U_q(\Hlie)$-submodule $W\subset V$ satisfying $V\simeq (W\otimes_{\CC}\CC(u))$, $\U_q(\hat{\Hlie})(W\otimes \CC[u^{\pm}])\subset (W\otimes \CC[u^{\pm}])$, and some additional technical properties. It
appears that the action of $\U_q(\hat{\Glie})$ has a nice regularity
property which is compatible with the Drinfeld coproduct and leads
to this rationality, that is to say that $\text{Mod}(\U_q'(\hat{\Glie}))$ is stable under a tensor product. The second point is the associativity : the
Drinfeld coproduct is coassociative, but the deformation parameter
$u$ breaks the symmetry and the $u$-deformed Drinfeld coproduct is
not coassociative in the usual sense. However there is a "twisted"
coassociative property which allows us to get the associativity of
a new monoidal structure on the category $\text{Mod}=\text{Mod}^0\oplus(\text{Mod}(\U_q'(\hat{\Glie}))\oplus\text{Mod}(\U_q'(\hat{\Glie}))\oplus\cdots)$ (here $\text{Mod}^0$ is an abelian semi-simple category with a unique simple object corresponding to a neutral object of the category).

\noindent Besides we prove that for a large class of quantum affinizations
(including quantum affine algebras and most quantum toroidal algebras) the subcategory of modules with a finite composition series is stable under the monoidal bifunctor.
This proof uses the specialization process described below and an investigation of the compatibility property between generalizations of Frenkel-Reshetikhin $q$-characters
and the monoidal bifunctor. We get in particular that all $l$-highest weight modules of the category $\text{Mod}(\U_q(\hat{\Glie}))$ have a finite composition series.

\noindent In order to go back to the usual category
$\text{Mod}(\U_q(\hat{\Glie}))$, we prove the existence of certain
forms for all cyclic modules ($\mathcal{A}$-forms) : such forms can
be specialized at $u=1$ and give cyclic $\U_q(\hat{\Glie})$-modules;
the proof uses the new rationality property. As a consequence, we
can construct an $l$-highest weight fusion $\U_q(\hat{\Glie})$-module
$V_1*_f V_2$ from two $l$-highest weight $\U_q(\hat{\Glie})$-modules
$V_1,V_2$. In particular a simple module of
$\text{Mod}(\U_q(\hat{\Glie}))$ is a quotient of a fusion module of
fundamental representations (analog to a result of Chari-Pressley
for the usual coproduct of quantum affine algebras). We establish a
cyclicity property of the fusion of any $l$-highest weight modules :
it is always an $l$-highest weight module. For quantum affine
algebras this property is very different from the properties of the
usual tensor product. It allows to control the "size" of the modules
obtained by $*_f$; in particular we can produce "big" integrable
$l$-highest weight modules whose existence was a priori not known.
We get moreover another bifunctor
$\otimes_d:\text{Modf}(\U_q(\hat{\Glie}))\times\text{Modf}(\U_q(\hat{\Glie}))\rightarrow\text{Modf}(\U_q(\hat{\Glie}))$
where $\text{Modf}(\U_q(\hat{\Glie}))$ is the category of finite
dimensional representations (in general $\otimes_d$ does not coincide with
$*_f$).

\noindent Let us describe an application. One important aspect of the representation theory of
quantum affine algebras is the existence of the Kirillov-Reshetikhin
modules whose characters are given by the certain explicit "fermionic" formulas. The
proof of the corresponding statement (the Kirillov-Reshetikhin
conjecture) can be obtained from certain $T$-systems originated from the theory of integrable systems. 
These $T$-systems can be interpreted in the form exact sequences involving tensor products of Kirillov-Reshetikhin modules; they were first established in \cite{Nad} for simply-laced quantum affine algebras (with the main result of \cite{Nab}) and then in \cite{her06} for all quantum affine algebras (with a different proof).
In the present paper we establish for a large class of general
quantum affinizations exact sequences involving fusion of
Kirillov-Reshetikhin modules : the usual tensor product can be
replaced in the general situation by the fusion constructed in this
paper. It is related to new generalized $T$-systems that we define
and establish. The existence of such new $T$-systems is a new
particular regularity property of the representation theory of
general quantum affinizations.

 We would like to mention other possible applications and new
developments : for example the quantum fusion tensor category should
lead to the description of the blocks of the category
$\text{Mod}(\U_q(\hat{\Glie}))$. Besides the generalized $T$-systems
established in this paper should lead to new fermionic formulas in
analogy to the case of quantum affine algebras and should have a
nice interpretation in terms of integrable systems. It would be
interesting to relate the quantum fusion to the Feigin-Loktev fusion
procedure defined for classical affine algebras \cite{fl}. We believe that the
deformation of the Drinfeld coproduct used is a particular case of a
more general framework involving some deformed Hopf structure on
$\ZZ$-graded algebras ("quantum Hopf vertex algebras" that we hope
to describe in another paper).  Eventually the tensor category
constructed in this paper should be related to certain new Hopf
algebras in the Tannaka-Krein reconstruction philosophy.

Let us describe the organization of this paper :

In section \ref{bck} we give the background for general quantum
affinizations and their representations and we recall the main results of
\cite{her04}. In section \ref{fuscat} we construct the "quantum
fusion" tensor category $\text{Mod}$ (theorem \ref{bif}) with the
help of the intermediate category $\text{Mod}(\U_q'(\hat{\Glie}))$
(definition \ref{ucat}) and establish the rationality property in $u$.
The case of fusion of Kirillov-Reshetikhin modules
of the quantum affine algebra $\U_q(\hat{sl_2})$ is studied
explicitly in subsection \ref{exkr}. In section \ref{saform}
we define and study the notion of $\mathcal{A}$-form (definition
\ref{aform}). In particular we prove in theorem \ref{aformcy}
that one can define $\mathcal{A}$-forms of cyclic modules of the
category $\text{Mod}(\U_q'(\hat{\Glie}))$. It gives rise to the two
following constructions : a fusion $l$-highest weight module from
two $l$-highest weight modules (definition \ref{defsf}), and in the
case of a quantum affine algebra a new bifunctor for the category of
finite dimensional representations (theorem \ref{defsd}). In section \ref{flsection} we prove for a large class of quantum
affinizations that the subcategory of finite length modules is
stable under the monoidal bifunctor (theorem \ref{prectens}) by
proving the compatibility with the $q$-characters (theorem
\ref{prod}) and with the combinatorial fusion product introduced in \cite{her04}; as a
by-product we get a representation theoretical interpretation of
this combinatorial fusion product (corollary \ref{inter}). In section \ref{tsyst} we describe several applications : first we prove
that a simple module of $\text{Mod}(\U_q(\hat{\Glie}))$ is a
quotient of the fusion of fundamental representations (proposition \ref{mainapp}); we establish a cyclicity property which allows to
control the "size" of this fusion module (theorem \ref{lhw} and
corollary \ref{conts}). Then we establish an exact sequence
involving fusion of Kirillov-Reshetikhin modules (theorem
\ref{exact}). It is related to new generalized $T$-systems (theorem
\ref{gentsyst}). In subsection \ref{furth} we address
complements on additional questions and further possible
applications.

 The main results of this paper were first announced in the "s\'eminaire quantique" in Strasbourg and at the conference
 "Representations of Kac-Moody Algebras and Combinatorics" in Banff respectively in January and March 2005.

{\bf Acknowledgments :} The author would like to thank V. Chari, E.
Frenkel, A. Moura, N. Reshetikhin, M. Rosso, O. Schiffmann and M.
Varagnolo for useful discussions. This paper was completed as the
author visited the university La Sapienza in Rome as a Liegrits
visitor; he would like to thank C. De Concini and C. Procesi for
their hospitality.

\section{Background}\label{bck} In this section we give backgrounds on general quantum affinizations and their representations (we also remind results of \cite{her04}). In the following for a formal variable $u$, $\CC[u^{\pm}],\CC(u),\CC[[u]],\CC((u))$ are the standard notations. 

\subsection{Cartan matrix}\label{datas} A generalized Cartan matrix is a matrix $C=(C_{i,j})_{1\leq i,j\leq n}$ satisfying \label{carmat}
$C_{i,j}\in\ZZ$, $C_{i,i}=2$, $(i\neq j\Rightarrow C_{i,j}\leq 0)$ and $(C_{i,j}=0\Leftrightarrow C_{j,i}=0)$. We denote
$I=\{1,...,n\}$ and $l=\text{rank}(C)$. For $i,j\in I$, we put $\delta_{i,j}=0$ if $i\neq j$, and $\delta_{i,j}=1$ if $i=j$.

\noindent In the following we suppose that $C$ is symmetrizable, that means that there is a matrix
$D=\text{diag}(r_1,...,r_n)$ ($r_i\in\NN^*$)\label{ri} such that $B=DC$\label{symcar} is symmetric. In particular if
$C$ is symmetric then it is symmetrizable with $D=I_n$. In the following $B$, $C$ and $D$ are fixed.

\noindent $q\in\CC^*$ is not a root of unity and is fixed. We put $q_i=q^{r_i}$ and for $l\in\ZZ$, we set
$[l]_q=\frac{q^l-q^{-l}}{q-q^{-1}}\in\ZZ[q^{\pm}]$. Let $C(q)$ be the quantized Cartan matrix defined by ($i\neq j\in
I$): $$C_{i,i}(q)=q_i+q_i^{-1}\text{ , }C_{i,j}(q)=[C_{i,j}]_q.$$
We denote by $D(q)$ the diagonal matrix such that for
$i,j\in I$, $D_{i,j}(q)=\delta_{i,j}[r_i]_q$. If $\text{det}(C(q))\neq 0$, we denote by $\tilde{C}(q)$ the inverse
matrix of $C(z)$.

\noindent We fix a realization $(\Hlie, \Pi, \Pi^{\vee})$ of $C$ (see \cite{kac}): $\Hlie$ is a $2n-l$
dimensional $\QQ$-vector space, $\Pi=\{\alpha_1,...,\alpha_n\}\subset \Hlie^*$ (set of the simple roots), $\Pi^{\vee}=\{\alpha_1^{\vee},...,\alpha_n^{\vee}\}\subset \Hlie$ (set of simple coroots) and for $1\leq i,j\leq n$: $\alpha_j(\alpha_i^{\vee})=C_{i,j}$. Denote by $\Lambda_1,...,\Lambda_n\in\Hlie^*$ (resp. the
$\Lambda_1^{\vee},...,\Lambda_n^{\vee}\in\Hlie$) the fundamental weights (resp. coweights). By definition $\alpha_i(\Lambda_j^{\vee})=\Lambda_i(\alpha_j^{\vee})=\delta_{i,j}$.

\noindent Consider a symmetric bilinear form $(,):\Hlie^*\times \Hlie^*\rightarrow \QQ$ such that for $i\in I$, $v\in\Hlie^*$ : $(\alpha_i,v)=v(r_i\alpha_i^{\vee})$. It is non degenerate and gives an isomorphism $\nu:\Hlie^*\rightarrow \Hlie$. In particular for $i\in I$ we have $\nu(\alpha_i)=r_i\alpha_i^{\vee}$, and for $\lambda,\mu\in\Hlie^*$ $\lambda(\nu(\mu))=\mu(\nu(\lambda))$.

\noindent Denote $P=\{\lambda \in\Hlie^*|\forall i\in I, \lambda(\alpha_i^{\vee})\in\ZZ\}$ the set of integral weights and $P^+=\{\lambda \in P|\forall i\in I, \lambda(\alpha_i^{\vee})\geq 0\}$ the set of dominant weights. For example we have $\alpha_1,...,\alpha_n\in P$ and $\Lambda_1,...,\Lambda_n\in P^+$. Denote $Q={\bigoplus}_{i\in I}\ZZ \alpha_i\subset P$ the root lattice and $Q^+={\sum}_{i\in I}\NN \alpha_i\subset Q$. Define $h:Q^+\rightarrow\NN$ such that $h(l_1\alpha_1+...+l_n\alpha_n)=l_1+...+l_n$.

\noindent For $\lambda,\mu\in \Hlie^*$, write $\lambda \geq \mu$ if $\lambda-\mu\in Q^+$. For $\lambda\in \Hlie^*$, denote $\mathcal{S}(\lambda)=\{\mu\in\Hlie^*|\mu\leq \lambda\}$.

\subsection{Quantum Kac-Moody algebra}\label{qkma} In the following $\Glie$ is the symmetrizable Kac-Moody algebra associated with the data $(\Hlie, \Pi,\Pi^{\vee})$ fixed in section \ref{datas} (see \cite{kac}).

\begin{defi} The quantum Kac-Moody algebra $\U_q(\Glie)$ is the $\CC$-algebra with generators $k_h$ ($h\in \Hlie$), $x_i^{\pm}$ ($i\in I$) and relations:
$$k_hk_{h'}=k_{h+h'}\text{ , }k_0=1\text{ , }k_hx_j^{\pm}k_{-h}=q^{\pm \alpha_j(h)}x_j^{\pm},$$
$$x_i^+x_j^- - x_j^-x_i^+=\delta_{i,j}\frac{k_{r_i\alpha_i^{\vee}}-k_{-r_i\alpha_i^{\vee}}}{q_i-q_i^{-1}},$$
$${\sum}_{r=0... 1-C_{i,j}}(-1)^r\begin{bmatrix}1-C_{i,j}\\r\end{bmatrix}_{q_i}(x_i^{\pm})^{1-C_{i,j}-r}x_j^{\pm}(x_i^{\pm})^r=0 \text{ (for $i\neq j$)}.$$
\end{defi}

\noindent This algebra was introduced in \cite{Dri1, jim}. It is remarkable that one can define a Hopf algebra structure on $\U_q(\Glie)$. In this paper we use the standard coproduct given by 
$$\Delta(k_h)=k_h\otimes k_h\text{ , }\Delta(x_i^+)=x_i^+\otimes 1 + k_i^+\otimes x_i^+\text{ , }\Delta(x_i^-)=x_i^-\otimes k_i^- + 1\otimes x_i^-,$$
where we denote $k_i^{\pm 1}=k_{\pm r_i\alpha_i^{\vee}}$. In particular we have a tensor structure on the category of $\U_q(\Glie)$-modules.

\subsection{Quantum affinizations}\label{genaff} For $\Glie$ a semi-simple Lie algebra, the affine Kac-Moody algebra associated to $\Glie$ is usually denoted by $\hat{\Glie}$ (see \cite{kac}). That is why the quantum affine algebra is usually denoted by $\U_q(\hat{\Glie})$. We will use this notation for general quantum affinizations. In the following $\Glie$ is the symmetrizable Kac-Moody algebra fixed in section \ref{datas} with datas $C,B,D$. We will define $\U_q(\hat{\Glie})$ the affinization of $\U_q(\Glie)$ as an algebra with generators $x_{i,r}^{\pm}$ ($i\in I, r\in\ZZ$), $k_h$ ($h\in \Hlie$), $h_{i,m}$ ($i\in I, m\in\ZZ-\{0\}$). Let $z$ be a formal variable and $\delta(z)={\sum}_{r\in\ZZ}z^r$. The relations are given in terms of the currents
$$x_i^{\pm}(w)={\sum}_{r\in\ZZ}x_{i,r}^{\pm}w^r,$$
$$\phi_i^{\pm}(z)={\sum}_{m\geq 0}\phi_{i,\pm m}^{\pm}z^{\pm m}=k_{\pm r_i \alpha_i^{\vee}}\exp(\pm(q-q^{-1}){\sum}_{m'\geq 1}h_{i,\pm m'}z^{\pm m'}).$$
We also set $\phi_{i,m}^+=0$ for $m<0$ and $\phi_{i,m}^-=0$ for $m>0$. (In this paper we omit the central elements because they act trivially on the considered representations).

\begin{defi}\label{defiaffi} The quantum affinization $\U_q(\hat{\Glie})$ is the $\CC$-algebra with generators $x_{i,r}^{\pm}$ ($i\in I, r\in\ZZ$), $k_h$ ($h\in \Hlie$), $h_{i,m}$ ($i\in I, m\in\ZZ-\{0\}$) and relations ($h,h'\in\Hlie$, $i,j\in I$) :
\begin{equation}\label{afcartu}k_hk_{h'}=k_{h+h'}\text{ , }k_0=1\text{ ,
}k_h\phi_i^{\pm}(z)=\phi_i^{\pm}(z)k_h,\end{equation}
\begin{equation}\label{actcartc}k_hx_j^{\pm}(z)=q^{\pm \alpha_j
(h)}x_j^{\pm}(z)k_h,\end{equation}
\begin{equation}\label{actcartplus}\phi_i^+(z)x_j^{\pm}(w)=\frac{q^{\pm
B_{i,j}}w-z}{w-q^{\pm B_{i,j}}z}x_j^{\pm}(w)\phi_i^+(z),\end{equation}
\begin{equation}\label{actcartmoins}\phi_i^-(z)x_j^{\pm}(w)=\frac{q^{\pm B_{i,j}}w-z}{w-q^{\pm
B_{i,j}}z}x_j^{\pm}(w)\phi_i^-(z),\end{equation}
\begin{equation}\label{partact}x_i^+(z)x_j^-(w)-x_j^-(w)x_i^+(z)=\frac{\delta_{i,j}}{q_i-q_i^{-1}}(\delta(\frac{w}{z})\phi_i^+(w)-\delta(\frac{z}{w})\phi_i^-(z)),\end{equation}
\begin{equation}\label{plusmoinsc}(w-q^{\pm B_{i,j}}z)x_i^{\pm}(z)x_j^{\pm}(w)=(q^{\pm
B_{i,j}}w-z)x_j^{\pm}(w)x_i^{\pm}(z),\end{equation}
\begin{equation}\label{equadeuxc}{\sum}_{\pi\in
\Sigma_s}{\sum}_{k=0...s}(-1)^k\begin{bmatrix}s\\k\end{bmatrix}_{q_i}x_i^{\pm}(w_{\pi(1)})...x_i^{\pm}(w_{\pi(k)})x_j^{\pm}(z)x_i^{\pm}(w_{\pi(k+1)})...x_i^{\pm}(w_{\pi(s)})=0,\end{equation}
where the last relation holds for all $i\neq j$, $s=1-C_{i,j}$, all sequences of integers $r_1,...,r_s$. $\Sigma_s$ is the symmetric group on $s$ letters. In the equation (\ref{actcartplus}) (resp. equation (\ref{actcartmoins})), $\frac{q^{\pm B_{i,j}}w-z}{w-q^{\pm B_{i,j}}z}$ is viewed as a formal power series in $z/w$ expanded around $\infty$ (resp. $0$).\end{defi}

\noindent When $C$ is of finite type, $\U_q(\hat{\Glie})$ is called a
quantum affine algebra, and when $C$ is of affine type,
$\U_q(\hat{\Glie})$ is called a quantum toroidal algebra. There is
a huge amount of very interesting papers on quantum affine algebras (see the introduction for references). They have the very particular property to be also quantum Kac-Moody
algebras \cite{Dri2, bec}. The quantum toroidal algebras
are also of particular importance and are closely related via a Frobenius-Schur-Weyl duality \cite{var1} to double affine Hecke algebras (Cherednik's algebras, see the introduction for references).
In general a quantum affinization (for example a quantum toroidal
algebra) is not isomorphic to a quantum Kac-Moody algebra.

\noindent Note that formulas (\ref{actcartplus}),
(\ref{actcartmoins}) are equivalent to (see for example
\cite{her04}) :
\begin{equation}\label{actcartld}h_{i,m}x_{j,r}^{\pm}-x_{j,r}^{\pm}h_{i,m}=\pm \frac{1}{m}[mB_{i,j}]_qx_{j,m+r}^{\pm}.\end{equation}

\noindent Relations (\ref{equadeuxc}) are called affine quantum Serre relations.

\noindent There is an injective algebra morphism $\U_q(\Glie)\rightarrow \U_q(\hat{\Glie})$ defined by $k_h\mapsto k_h$ , $x_i^{\pm}\mapsto x_{i,0}^{\pm}$ for $h\in \Hlie, i\in I$. For $i\in I$, the subalgebra $\hat{\U}_i\subset \U_q(\hat{\Glie})$ generated by the $x_{i,r}^{\pm}, h_{i,m}, k_{p \alpha_i^{\vee}}$ ($r\in\ZZ, m\neq 0, p\in\QQ$) is isomorphic to $\U_{q_i}(\hat{sl_2})$.

\noindent Let $\U_q(\hat{\Glie})^{\pm}$ (resp. $\U_q(\hat{\Hlie})$) be the subalgebra of
$\U_q(\hat{\Glie})$ generated by the $x_{i,r}^{\pm}$ (resp. by the $k_h$ and the $h_{i,r}$). We have a triangular decomposition :

\begin{thm}\cite{bec, her04}\label{afftrian} The multiplication  $x^-\otimes h\otimes x^+\mapsto x^-hx^+$ defines an isomorphism of $\CC$-vector space $\U_q(\hat{\Glie})^-\otimes \U_q(\hat{\Glie})\otimes \U_q(\hat{\Glie})^+\simeq \U_q(\hat{\Glie})$.\end{thm}

\noindent The result was proved in \cite[Proposition 6.1]{bec} for quantum affine algebras (see also \cite[Proposition 12.2.2]{Cha2}) and the general case is treated in \cite[Theorem 3.2]{her04} (the result follows from \cite[Proposition 6.1]{bec} for Cartan matrices satisfying $(i\neq j\Rightarrow C_{i,j}C_{j,i}\leq 3)$; in general the technical point is to check the compatibility with affine quantum Serre relations (\ref{equadeuxc}). The proof \cite{her04} implies some combinatorial relations discovered by Jing \cite[Theorem 6.1]{jin}).

\subsection{The category $\text{Mod}(\U_q(\hat{\Glie}))$}\label{int} In this section we recall results on the category of integrable $\U_q(\hat{\Glie})$-modules which are in the category $\mathcal{O}$ defined below.

\noindent For $V$ a $\U_q(\Hlie)$-module and $\omega\in \Hlie^*$ we denote by $V_{\omega}$ the weight space of weight $\omega$
$$V_{\omega}=\{v\in V|\forall h\in \Hlie, k_h.v=q^{\omega(h)}v\}.$$

\noindent We say that $V$ is $\U_q(\Hlie)$-diagonalizable if $V={\bigoplus}_{\omega\in \Hlie^*}V_{\omega}$. We have the following standard definitions (see \cite{kac} for the classical definition) :

\begin{defi} A $\U_q(\Hlie)$-module $V$ is said to be integrable if $V$ is $\U_q(\Hlie)$-diagonalizable and for all $\omega\in \Hlie^*$, $V_{\omega}$ is finite dimensional, and for $i\in I$ there is $R\geq 0$ such that $(r\geq R \Rightarrow V_{\omega\pm r\alpha_i}=\{0\})$. A $\U_q(\hat{\Glie})$-module (resp. a $\U_q(\Glie)$-module) $V$ is said to be integrable if $V$ is integrable as a $\U_q(\Hlie)$-module.\end{defi}

\noindent An integrable representation is not necessarily finite dimensional.

\noindent We have the natural analog of the classical category $\mathcal{O}$ of Kac-Moody algebras \cite{kac} :

\begin{defi} A $\U_q(\Hlie)$-module $V$ is said to be in the category $\mathcal{O}(\U_q(\Hlie))$ if

i) $V$ is $\U_q(\Hlie)$-diagonalizable,

ii) for all $\omega\in \Hlie^*$, $\text{dim}(V_{\omega})<\infty$,

iii) there is a finite number of element $\lambda_1,...,\lambda_s\in \Hlie^*$ such that the weights of $V$ are in ${\bigcup}_{j=1,\cdots,s}\mathcal{S}(\lambda_j)$.\end{defi}

\begin{defi} For $\mathfrak{m}$ respectively equal to $\hat{\Hlie}$, $\Glie$, $\hat{\Glie}$, we denote by $\text{Mod}(\U_q(\mathfrak{m}))$ the category of $\U_q(\mathfrak{m})$-modules which are integrable and in the category $\mathcal{O}(\U_q(\Hlie))$ as a $\U_q(\Hlie)$-module.\end{defi}

\noindent The category $\text{Mod}(\U_q(\Glie))$ (in a slightly more general form) is considered in \cite{lu}. Definitions for $\U_q(\hat{\Glie})$-modules are given in \cite[section 1.2]{Naams}.

\noindent For $\lambda\in P$, denote by $L(\lambda)$ the
simple highest weight $\U_q(\Glie)$-module of highest weight
$\lambda$. We have (see for example \cite{ro1, lu}) :

\begin{thm}\label{usual} $L(\lambda)\in\text{Mod}(\U_q(\Glie))$ if and only if $\lambda\in P^+$.\end{thm}

The category $\text{Mod}(\U_q(\Glie))$ is semi-simple. But
$\text{Mod}(\U_q(\hat{\Glie}))$ is not semi-simple. 

By generalizing the approach of finite dimensional representations of quantum affine
algebras developed in \cite{Cha2}, and definitions in \cite{Naams, mi}, we defined in \cite{her04} for
general quantum affinization $\U_q(\hat{\Glie})$ :

\begin{defi}\label{aint} An $l$-weight is a couple $(\lambda,\Psi)$ where $\lambda\in P$ and $\Psi=(\Psi_{i,\pm m}^{\pm})_{i\in I, m\geq 0}$, such that $\Psi_{i,\pm m}^{\pm}\in\CC^*$ and $\Psi_{i,0}^{\pm}=q_i^{\pm \lambda(\alpha_i^{\vee})}$ for $i\in I$ and $m\geq 0$.

\noindent A $\U_q(\hat{\Glie})$-module $V$ is said to be of $l$-highest weight $(\lambda,\Psi)$ if there is $v\in V_{\lambda}$ ($l$-highest weight vector) such that $V=\U_q(\hat{\Glie}).v$ and $x_{i,r}^+.v=0$, $\phi_{i,\pm m}^{\pm}.v=\Psi_{i,\pm m}^{\pm}v$ for $i\in I$, $r\in\ZZ$ and $m\geq 0$.

\noindent An $l$-weight $(\lambda,\Psi)$ is said to be dominant if there is an $n$-tuplet of (Drinfeld) polynomials $(P_i(u))_{i\in I}\in (\CC[u])^n$ satisfying for $i\in I$, $P_i(0)=1$ and the relation in $\CC[[z]]$ (resp. in $\CC[[z^{-1}]]$):
\begin{equation}\label{defipsi}
{\sum}_{m\geq 0} \Psi_{i,\pm m}^{\pm} z^{\pm m}=q_i^{\text{deg}(P_i)}\frac{P_i(zq_i^{-1})}{P_i(zq_i)}.
\end{equation}
\end{defi}

\noindent The element $\Psi$ in an $l$-weight $(\lambda,\Psi)$ is called a pseudo-weight.

\noindent Theorem \ref{afftrian} implies that for $(\lambda,\Psi)$ an $l$-weight there is a unique simple module $L(\lambda,\Psi)$ of $l$-highest weight $(\lambda,\Psi)$. Moreover :

\begin{thm}\cite{Cha2, mi, Naams, her04}\label{intsimp} We have $L(\lambda,\Psi)\in \text{Mod}(\U_q(\hat{\Glie}))$ if and only if $(\lambda,\Psi)$ is dominant.\end{thm}

\noindent This result was proved in \cite[Theorem 12.2.6]{Cha2} for quantum affine
algebras (moreover in this case these integrable representations are
finite dimensional), in \cite[Theorem 1]{mi} for quantum toroidal algebras of
type $A$, in \cite[Proposition 1.2.15]{Naams} for general simply-laced quantum affinizations and
in \cite[Theorem 4.9]{her04} for general quantum affinizations.

\noindent {\it Examples.} For $i\in I, a\in\CC^*, r\geq 1$, consider the $l$-weight $(\Lambda_i,\Psi_{i,a,r})$ where $\Psi_{i,a,r}$ is given by equation (\ref{defipsi}) with the $n$-tuplet $P_i(u)=(1-ua)(1-uaq_i^2)\cdots(1-uaq_i^{2(r-1)})$ and $P_j(u)=1$ for $j\neq i$. The $\U_q(\hat{\Glie})$-module $W_{r,a}^{(i)}=L(\Lambda_i,\Psi_{i,a,r})\in \text{Mod}(\U_q(\hat{\Glie}))$ is called a Kirillov-Reshetikhin module. 

\noindent In the particular case $r=1$, the $\U_q(\hat{\Glie})$-module $L_{i,a}=W_{1,a}^{(i)}$ is called a fundamental representation. 

\noindent For $\lambda\in P$ satisfying ($\forall j\in I, \lambda(\alpha_j^{\vee})=0$), consider the $l$-weight $(\lambda,\Psi_0)$ where $\Psi_0$ is given by equation (\ref{defipsi}) with the $n$-tuplet $P_j(u)=1$ for all $j\in I$. The $\U_q(\hat{\Glie})$-module $L_{\lambda}=L(\lambda,\Psi_0)\in\text{Mod}(\U_q(\hat{\Glie}))$ is also called a fundamental representation in this paper. 

\section{Construction of the quantum fusion tensor category $\text{Mod}$}\label{fuscat} 
In general a quantum affinization $\U_q(\hat{\Glie})$ has no Hopf algebra structure (except in the
case of quantum affine algebras where we have the coproduct of the
Kac-Moody realization). However Drinfeld (unpublished note, see also
\cite{di, df}) defined for $\U_q(\hat{sl_n})$ a map which behaves like
a new coproduct adapted to the affinization realization, but it is
defined in a completion and can not directly be used to define
tensor product of representations. In this paper we construct a
corresponding tensor category (section \ref{fuscat}) by deforming the Drinfeld coproduct, see how to go
back to the category $\text{Mod}(\U_q(\hat{\Glie}))$ by
specialization (section \ref{saform}), and give applications of
these results (section \ref{tsyst}). We also prove in section
\ref{flsection} that it induces a tensor category structure on the
subcategory of modules of finite length. For quantum affine algebras
this process is new and different from the usual tensor product (we
prove that it has different properties), and for general quantum
affinizations (as general quantum toroidal algebras) there was no
process to construct fusion of modules so far. For this construction, the
first technical point that we solve is a rationality problem for the
dependence in the deformation parameter $u$ of the Drinfeld coproduct 
(lemma \ref{proddef}). The second point is the
associativity : the $u$-deformed Drinfeld coproduct has a "twisted"
coassociative property which allows us to get the associativity of
the new monoidal structure (lemma \ref{coass}).

\noindent In this section we define the category $\text{Mod}$
(definition \ref{ucat}) and the tensor structure (theorem
\ref{bif}). The case of fusion of Kirillov-Reshetikhin modules of
the quantum affine algebra $\U_q(\hat{sl_2})$ is studied explicitly
in subsection \ref{exkr}.

\subsection{Deformation of the Drinfeld coproduct}\label{udelta} In this section we study a $u$-deformation of the Drinfeld coproduct, in particular a "twisted" coassociativity property (lemma \ref{coass}).

\noindent Let $\U_q'(\hat{\Glie})=\U_q(\hat{\Glie})\otimes \CC(u)$ and
$\U_q'(\hat{\Glie})\hat{\otimes}\U_q'(\hat{\Glie})=(\U_q(\hat{\Glie})\otimes_{\CC}\U_q(\hat{\Glie}))((u))$ be the $u$-topological completion of $\U_q'(\hat{\Glie})\otimes_{\CC(u)}\U_q'(\hat{\Glie})$. Let $\tilde{\U}_q(\hat{\Glie})$ be the algebra defined as $\U_q(\hat{\Glie})$ without affine quantum Serre relations (\ref{equadeuxc}). We also define $\tilde{\U}_q'(\hat{\Glie})=\tilde{\U}_q(\hat{\Glie})\otimes \CC(u)$ and $\tilde{\U}_q'(\hat{\Glie})\hat{\otimes}\tilde{\U}_q'(\hat{\Glie})=(\tilde{\U}_q(\hat{\Glie})\otimes_{\CC}\tilde{\U}_q(\hat{\Glie}))((u))$. 

In \cite[Proposition 6.3]{her04} we introduced a $u$-deformation of the Drinfeld coproduct :

\begin{prop}\label{coprod} There is a unique morphism of $\CC(u)$-algebra $\Delta_u: \tilde{\U}_q'(\hat{\Glie})\rightarrow \tilde{\U}_q'(\hat{\Glie})\hat{\otimes}\tilde{\U}_q'(\hat{\Glie})$ such that for $i\in I$, $r\in\ZZ$, $m\geq 0$, $h\in\Hlie$:
\begin{equation}\label{deltaup}\Delta_u(x_{i,r}^+)= x_{i,r}^+\otimes 1 +{\sum}_{s\geq 0} u^{r+l} (\phi_{i,-s}^-\otimes x_{i,r+s}^+),\end{equation}
\begin{equation}\label{deltaum}\Delta_u(x_{i,r}^-)=u^r (1\otimes x_{i,r}^-) + {\sum}_{s\geq 0} u^s (x_{i,r-s}^-\otimes \phi_{i,s}^+),\end{equation}
\begin{equation}\label{deltauh}\Delta_u(\phi_{i,\pm m}^{\pm})= {\sum}_{0\leq s\leq m}u^{\pm s} (\phi_{i,\pm (m-s)}^{\pm}\otimes\phi_{i,\pm s}^{\pm}),\end{equation}
\begin{equation}\label{deltauh2}\Delta_u(k_h)=k_h\otimes k_h.\end{equation}
\end{prop}

\noindent If $(i\neq j\Rightarrow C_{i,j}C_{j,i}\leq 3)$, then $\Delta_1$ is compatible with affine quantum Serre relations (\ref{equadeuxc}) (see \cite{di, e, gr}). By definition of $\Delta_u$, this can be easily generalized to :

\begin{cor}\label{mtrois} If $(i\neq j\Rightarrow C_{i,j}C_{j,i}\leq 3)$, then the map $\Delta_u$ induces a morphism $\Delta_u:
\U_q'(\hat{\Glie})\rightarrow \U_q'(\hat{\Glie})\hat{\otimes}\U_q'(\hat{\Glie})$.\end{cor}

\noindent This case includes quantum affine algebras and most quantum toroidal algebras (except $A_1^{(1)}, A_{2l}^{(2)}$, $l\geq 0$).

\begin{rem}\label{toref} Note that the simple $l$-highest weight modules of $\tilde{\U}_q(\hat{\Glie})$ are the same as for $\U_q(\hat{\Glie})$ (via the projection $\tilde{\U}_q(\hat{\Glie})\rightarrow\U_q(\hat{\Glie})$). From corollary \ref{mtrois}, if ($i\neq j\Rightarrow C_{i,j}C_{j,i}\leq 3$), the results in this paper can be indifferently stated and proved for $\U_q(\hat{\Glie})$ or $\tilde{\U}(\hat{\Glie})$. So in some sections, we use the following notation :

1) If $(i\neq j\Rightarrow C_{i,j}C_{j,i}\leq 3$), $\U_q(\hat{\Glie})$ (resp. $\U_q(\Glie), \U_q'(\hat{\Glie})$) means the algebra with affine quantum Serre relations (\ref{equadeuxc}).

2) Otherwise, $\U_q(\hat{\Glie})$ (resp. $\U_q(\Glie), \U_q'(\hat{\Glie})$) means the algebra without affine quantum Serre relations (\ref{equadeuxc}).

\noindent This point will be explicitly reminded in the rest of this paper by references to this remark.\end{rem}

\noindent Observe that for $i\in I$, the subalgebra $\hat{\U}_i'=\hat{\U}_i\otimes_{\CC}\CC(u)\subset \U_q'(\hat{\Glie})$ is a "Hopf subalgebra" of $\U_q'(\hat{\Glie})$ for $\Delta_u$, that is to say $\Delta_u(\hat{\U}_i')\subset \hat{\U}_i'\hat{\otimes}\hat{\U}_i'$.

The usual coassociativy property of a coproduct $\Delta$ is $(\text{Id}\otimes\Delta)\circ\Delta=(\Delta\otimes\text{Id})\circ\Delta$. This relation is not satisfied by $\Delta_u$. However, we have a "twisted" coassociativity property satisfied by $\Delta_u$ and obtained by replacing the parameter $u$ by some power of $u$ (at the "limit" $u=1$ we recover the usual coassociativity property). This relation will be crucial in the construction of the quantum fusion tensor category (note that we do not use the quasi-Hopf algebras point of view) :

\begin{lem}\label{coass} Let $\U_q(\hat{\Glie})$ as in remark \ref{toref}. Let $r,r'\geq 1$. As algebra morphisms $\U_q(\hat{\Glie})\rightarrow(\U_q(\hat{\Glie})\otimes\U_q(\hat{\Glie})\otimes\U_q(\hat{\Glie}))((u))$, the two following maps are equal :
$$(\text{Id}\otimes\Delta_{u^{r'}})\circ\Delta_{u^r}=(\Delta_{u^r}\otimes\text{Id})\circ\Delta_{u^{r+r'}}.$$
\end{lem}

\demo It suffices to check the equality on the generators. Because of $u$, the images of both applications are of the form 
$${\sum} u^{r(\text{deg}(g_{(2)})+\text{deg}(g_{(3)}))+r'\text{deg}(g_{(3)})} g_{(1)}\otimes g_{(2)}\otimes g_{(3)}$$
where $g_{(1)} , g_{(2)} , g_{(3)}$ are homogeneous. We just give the results of the maps on generators and leave it to the reader to check it. For $x_{i,p}^+$, $x_{i,p}^-$, $\phi_{i,\pm m}^{\pm}$, $k_h$, respectively the two maps give :
$$x_{i,p}^+\otimes 1\otimes 1+{\sum}_{s\geq 0} u^{r(p+s)} (\phi_{i,-s}^-\otimes x_{i,p+s}^+\otimes 1)+{\sum}_{s,s'\geq 0}u^{r(p+s)+r'(p+s+s')}(\phi_{i,-s}^-\otimes \phi_{i,-s'}^-\otimes x_{i,p+s+s'}^+),$$
$$u^{(r+r')p} (1\otimes 1 \otimes x_{i,p}^-) + {\sum}_{s\geq 0} u^{rp+r's} (1\otimes x_{i,p-s}^-\otimes \phi_{i,s}^+)+{\sum}_{s,s'\geq 0}u^{r(s+s')+s'r'}(x_{i,p-s-s'}^-\otimes \phi_{i,s}^+\otimes \phi_{i,s'}^+),$$
$${\sum}_{\{s_1,s_2,s_3\geq 0|s_1+s_2+s_3=m\}}u^{\pm ((s_2+s_3)r+s_3'r')} (\phi_{i,\pm s_1}^{\pm}\otimes\phi_{i,\pm s_2}^{\pm}\otimes\phi_{i,\pm s_3}^{\pm}),$$
$$k_h\otimes k_h\otimes k_h.$$\qed

\subsection{The category $\text{Mod}(\U_q'(\hat{\Glie}))$}

\noindent In this section we introduce and study the category $\text{Mod}(\U_q'(\hat{\Glie}))$ which is "a building graded part" of the quantum fusion tensor category.

\noindent First let us give a $u$-deformation of the notion of $l$-weight :

\begin{defi} An $(l,u)$-weight is a couple $(\lambda,\Psi(u))$ where $\lambda\in \Hlie^*$, $\Psi(u)=(\Psi_{i,\pm m}^{\pm}(u))_{i\in I, m\geq 0}$, $\Psi_{i,\pm m}^{\pm}(u)\in\CC[u^{\pm 1}]$ satisfying $\Psi_{i,0}^{\pm}(u)=q_i^{\pm \lambda(\alpha_i^{\vee})}$.\end{defi}

\begin{defi}\label{ucat} Let $\text{Mod}(\U_q'(\hat{\Glie}))$ be the category of $\U_q'(\hat{\Glie})$-modules $V$ such that there is a $\CC$-vector subspace $W\subset V$ satisfying :

i) $V\simeq W\otimes_{\CC}\CC(u)$,

ii) $W$ is stable under the action of $\U_q(\Hlie)$ and is an object of $\text{Mod}(\U_q(\Hlie))$,

iii) the image of $W\otimes \CC[u^{\pm}]$ in $V$ by the morphism of i) is stable under the action of $\U_q(\hat{\Hlie})$,

iv) for all $\omega\in\Hlie^*$, all $i\in I$, all $r\in\ZZ$, there
are a finite number of $\CC$-linear operators
$f_{k,k',\lambda}^{\pm} : W_{\omega}\rightarrow V$ ($k\geq 0, k'\geq
0, \lambda\in\CC^*$) such that for all $m\geq 0$ and all $v\in
W_{\omega}$ :
$$x_{i,r\pm m}^{\pm}(v)={\sum}_{k\geq 0,k'\geq 0,\lambda\in\CC^*} \lambda^m u^{km} m^{k'} f_{k,k',\lambda}^{\pm}(v),$$

v) we have a decomposition 
$$V={\bigoplus}_{(\lambda,\Psi(u))\text{ $(l,u)$-weight}}V_{(\lambda,\Psi(u))}$$
$$\text{where }V_{(\lambda,\Psi(u))}=\{v\in V_{\lambda}|\exists p\geq 0 , \forall i\in I, m\geq 0, (\phi_{i,\pm m}^{\pm}-\Psi_{i,\pm m}^{\pm}(u))^p.v=0\}.$$\end{defi}

\noindent We have a first result :

\begin{lem}\label{phim} Let $V\in\text{Mod}(\U_q'(\hat{\Glie}))$. For all $\omega\in\Hlie^*$, all $i\in I$, there are a finite number of $\CC$-linear operators
$g_{k,k',\lambda}^{\pm} :W_{\omega}\rightarrow V$ ($k\geq 0, k'\geq
0, \lambda\in\CC^*$) such that for all $m\geq 1$ and all $v\in
W_{\omega}$ :
$$\phi_{i,r\pm m}^{\pm}(v)={\sum}_{k\geq 0,k'\geq 0,\lambda\in\CC^*} \lambda^m u^{km} m^{k'} g_{k,k',\lambda}^{\pm}(v).$$\end{lem}

\demo Formula (\ref{partact}) gives for $m\geq 1$ :
$$\phi_{i,-m}^{-}=(q_i^{-1}-q_i)(x_{i,0}^+x_{i,-m}^- -x_{i,-m}^-x_{i,0}^+)\text{ and }\phi_{i,m}^+=(q_i-q_i^{-1})(x_{i,m}^+x_{i,0}^- -x_{i,0}^-x_{i,m}^+)$$
and so property iv) of definition \ref{ucat} gives the result.\qed

Now we give a typical example of an object of $\text{Mod}(\U_q'(\hat{\Glie}))$. For $V$ a module in $\text{Mod}(\U_q(\hat{\Glie}))$, consider $i(V)$ the $\U_q'(\hat{\Glie})$-module obtained by extension : $i(V)=V\otimes\CC(u)$. 

\begin{prop}\label{inclu} $i$ defines an injective faithful functor $i:\text{Mod}(\U_q(\hat{\Glie}))\rightarrow\text{Mod}(\U_q'(\hat{\Glie}))$.\end{prop}

\noindent In particular $\text{Mod}(\U_q(\hat{\Glie}))$ can be viewed as a subcategory of $\text{Mod}(\U_q'(\hat{\Glie}))$.

\demo Let us prove that the functor is well-defined. Consider $V\in \text{Mod}(\U_q(\hat{\Glie}))$ and let us prove that $i(V)\in\text{Mod}(\U_q'(\hat{\Glie}))$ where we choose $W=V$ for the $\CC$-vector space of definition
\ref{ucat}. The properties $i), ii), iii)$ are clear. The property $v)$ is clear because the base field for $V$ is $\CC$. Let us prove the property iv) : let $i\in I$,
$\omega\in\Hlie^*$ and let us consider the operators $x_{i,r}^+$, $r\in\ZZ$ (the proof for the $x_{i,r}^-$ is analog). As $W$ is integrable, there is a finite
dimensional $\hat{\U}_i$-submodule $W'$ such that $W_{\lambda}\subset W'\subset W$. Consider $\rho:\hat{\U}_i\rightarrow \text{End}(W')$ the action. Consider the linear map
$\Phi:\text{End}(W')\rightarrow\text{End}(W')$ defined by $\Phi(p)=\frac{1}{[2r_i]_q}(\rho(h_{i,1})p-p\rho(h_{i,1}))$. As $\text{End}(W')$ is finite dimensional, we can consider the
Jordan decomposition $\Phi=\Phi_1+\Phi_2$ where $\Phi_1\circ\Phi_2=\Phi_2\circ\Phi_1$, $\Phi_1$ is diagonalizable and $\Phi_2^{S}=0$ for some $S\geq 1$. It
follows from formula (\ref{actcartld}) that for $m\geq 0$, $\rho(x_{i,r+m}^+)=\Phi^m(\rho(x_{i,r}^+))$, and so
$$\rho(x_{i,r+m}^+)={\sum}_{s=0,\cdots,S}\begin{bmatrix}m\\s\end{bmatrix}(\Phi_2^s\Phi_1^{m-s})(\rho(x_{i,r}^+))={\sum}_{\lambda\in\CC}{\sum}_{s=0,\cdots,S}\begin{bmatrix}m\\s\end{bmatrix}\lambda^{m-s}(\Phi_2^s)(d_{\lambda})$$
where $\rho(x_{i,r}^+)={\sum}_{\lambda\in\CC}d_{\lambda}$ and for $\lambda\in\CC$, $\Phi_1(d_{\lambda})=\lambda d_{\lambda}$. In particular we get property $iv)$ of definition
\ref{ucat} because for a fixed $s\geq 0$, $\begin{bmatrix}m\\s\end{bmatrix}$ is a polynomial in $m$. For $m\leq 0$, we replace $\Phi$ by
$\Phi'(p)=\frac{1}{[2r_i]_q}(\rho(h_{i,-1})p-p\rho(h_{i,-1}))$. So $i$ is well-defined. It is clearly injective and faithful. \qed

\subsection{Construction of the tensor structure}\label{fusumod} The aim of this section is to define a tensor category $\text{Mod}$ by using $\Delta_u$ and $\text{Mod}(\U_q'(\hat{\Glie}))$. It is the main tool used in this paper; the stability inside the field $\CC(u)$ is
one of the crucial points that make the category $\text{Mod}$ useful
for the purposes of the present paper. 

\subsubsection{The category $\text{Mod}$} Let $\text{Mod}^0$ be the full abelian semi-simple subcategory of $\text{Mod}(\U_q'(\hat{\Glie}))$-modules with a unique simple object $i(L_0)$ (we recall that $L_0=L(0,\Psi)$ where $\Psi$ is defined with Drinfeld polynomials equal to $1$).

\begin{defi} We denote by $\text{Mod}$ the direct sum of categories 
$$\text{Mod}=\text{Mod}^0\oplus (\text{Mod}(\U_q'(\hat{\Glie}))\oplus \text{Mod}(\U_q'(\hat{\Glie}))\oplus \cdots).$$ 
For $r\geq 1$, the $r$-th summand in the second sum is denoted by $\text{Mod}^r$ : $\text{Mod}=\bigoplus_{r\geq 0}\text{Mod}^r$.\end{defi}

\noindent Note that with the identification $\text{Mod}(\U_q'(\hat{\Glie}))\simeq \text{Mod}^1$, we can also consider that $i$ is an injective faithful functor $\text{Mod}(\U_q(\hat{\Glie}))\rightarrow \text{Mod}$. So one can view $\text{Mod}(\U_q(\hat{\Glie}))$ as a subcategory of $\text{Mod}$.

\subsubsection{Rationality property} We use the notations of remark \ref{toref}. Fix $r\geq 1$. Let $V_1,V_2\in\text{Mod}(\U_q'(\hat{\Glie}))$. One defines an action of $\U_q'(\hat{\Glie})$ on $(V_1\otimes_{\CC(u)} V_2)((u))$ by the following formula ($g\in\U_q'(\hat{\Glie})$, $v_1\in V_1, v_2\in V_2$) :
\begin{equation}\label{defit} g.(v_1\otimes v_2)=\Delta_{u^r}(g)(v_1\otimes v_2).\end{equation}
(A priori this formula only makes sense for $r\geq 1$, but we will see bellow that in some particular situations we can also use it for $r=0$).

We have the following "rationality property" of the action given by formula (\ref{defit}), which is the crucial point for the construction of the tensor structure.

\begin{lem}\label{proddef} The subspace $(V_1\otimes_{\CC(u)} V_2)\subset(V_1\otimes_{\CC(u)} V_2)((u))$ is stable under the action of $\U_q'(\hat{\Glie})$ defined by formula (\ref{defit}). The induced $\U_q'(\hat{\Glie})$-module structure on $(V_1\otimes_{\CC(u)} V_2)$ is an object of $\text{Mod}(\U_q'(\hat{\Glie}))$.\end{lem}

\demo From definition \ref{ucat} we have subspaces $W_1\subset V_1, W_2\subset V_2$. We choose $W=W_1\otimes_{\CC} W_2 \subset (V_1\otimes_{\CC(u)} V_2)$ and we prove that the properties of definition \ref{ucat} are satisfied. properties
i), ii) are clear. Properties iii), v) follow directly from above formulas (\ref{deltauh}), (\ref{deltauh2}). Let us prove property iv). Let $\lambda,\mu\in\Hlie^*$, $p\in\ZZ$ and suppose that for $m\geq 0, v_1\in(W_1)_{\lambda}$, $x_{i,p\pm m}^{\pm}(v_1)=u^{mk_{1,\pm}}\lambda_{1,\pm}^m m^{k_{1,\pm}'} f_1^{\pm}(v_1)$, for
$m\geq 0, v_2\in(W_2)_{\mu}$, $x_{i,p\pm m}^{\pm}(v_2)=u^{mk_{2,\pm}}\lambda_{2,\pm}^m m^{k_{2,\pm}'} f_2^{\pm}(v_2)$. It follows
from lemma \ref{phim}, that we can suppose that for $t\geq 1$, $\phi_{i,-t}^-(v_1)=u^{k_3t}\lambda_3^t t^{k_3'} f_3(v)$. Let $v_1\in
(W_1)_{\lambda}, v_2\in (W_2)_{\mu}$ and $v=v_1\otimes v_2$. Formulas (\ref{deltaup}) and (\ref{defit}) give $x_{i,p+m}^+(v)=A+B$ where 
$$A=x_{i,p+m}^+(v_1)\otimes v_2+u^{r(p+m)} q_i^{-\lambda(\alpha_i^{\vee})}v_1\otimes x_{i,p+m}^+(v_2)$$ $$=u^{mk_{1,+}}\lambda_{1,+}^m
m^{k_{1,+}'} f_1^+(v_1)\otimes v_2 +u^{r(p+m)} q_i^{-\lambda (\alpha_i^{\vee})}v_1\otimes u^{mk_{2,+}}\lambda_{2,+}^m m^{k_{2,+}'} f_2^+(v_2),$$ 
and 
$$B={\sum}_{t\geq 1} u^{r(p+m+t)} (\phi_{i,-t}^-(v_1)\otimes x_{i,p+m+t}^+(v_2))$$ $$={\sum}_{t\geq 1} u^{r(p+m+t)} u^{k_3t}\lambda_3^t t^{k_3'}
f_3(v_1)\otimes u^{k_{2,+}(m+t)}\lambda_{2,+}^{m+t} (m+t)^{k_{2,+}'} f_2^{+}(v_2)$$
$$={\sum}_{s=0,\cdots,k_{2,+}'}\begin{bmatrix}k_{2,+}'\\s\end{bmatrix}\lambda_{2,+}^mu^{r(p+m)}u^{mk_{2,+}}m^{k_{2,+}'-s}\mathcal{R}_s(u) f_3(v_1)\otimes
f_2^{+}(v_2),$$ 
where 
$$\mathcal{R}_s(u)={\sum}_{t\geq 1}
u^{rt} u^{k_3t}\lambda_3^t u^{tk_{2,+}}\lambda_{2,+}^tt^{s+k_3'}.$$
As $\mathcal{R}_s(u)\in\CC(u)$, $x_{i,p+m}^+(v)$ makes sense in
$V_1\otimes_{\CC(u)} V_2$. Moreover $\mathcal{R}_s(u)$ does not
depend of $m$, and so it follows from the above formulas for $A$ and $B$ that property iv) of definition \ref{ucat} is satisfied (the study is analog for $x_{i,m}^-$). So the $\U_q'(\hat{\Glie})$-module $V_1\otimes_{\CC(u)} V_2$ is an object of $\text{Mod}(\U_q'(\hat{\Glie}))$.\qed

\subsubsection{Definition of the tensor structure} We use the notations of remark \ref{toref}. First let us study how formula (\ref{defit}) behaves with the representation $i(L_0)$.

\begin{lem}\label{extend} Let $r\geq 1$ and $V\in\text{Mod}^r$. Then formula (\ref{defit}) defines a structure of $\U_q'(\hat{\Glie})$-module on $V\otimes_{\CC(u)} i(L_0)$ (resp. on $i(L_0)\otimes_{\CC(u)} V$) which is isomorphic to $V$.\end{lem}

\demo For $i\in I, s\in\ZZ$ we have $x_{i,s}^{\pm}.i(L_0)=\{0\}$ and for $i\in I, s>0$ we have $\phi_{i,\pm s}^{\pm}.i(L_0)=\{0\}$. For $i\in I$, the action of $\phi_{i,0}^{\pm}$ on $i(L_0)$ is the identity. So for $V\otimes_{\CC(u)} i(L_0)$, the action given by formula (\ref{defit}) is $g\otimes 1$ which makes sense for $g\in\U_q'(\hat{\Glie})$ (direct computation on generators). For $i(L_0)\otimes_{\CC(u)} V$ the action given by formula (\ref{defit}) is $1\otimes g$ which makes sense for $g\in\U_q'(\hat{\Glie})$.\qed

Fix $r,r'\geq 0$. Let $V_1\in\text{Mod}^r$ and $V_2\in\text{Mod}^{r'}$. If $r,r'\geq 1$, one defines an action of $\U_q'(\hat{\Glie})$ on $V_1\otimes_{\CC(u)} V_2$ by formula (\ref{defit}). From lemma \ref{proddef} we get an object in $\text{Mod}(\U_q'(\hat{\Glie}))$ (this can be extended to the cases $r=0$ or $r'=0$ by lemma \ref{extend}). We consider this tensor product as an object in the $(r+r')$-th summand of $\text{Mod}$. So we have defined a bifunctor $\otimes_f:\text{Mod}\times\text{Mod}\rightarrow\text{Mod}$.

\begin{thm}\label{bif}  The bifunctor $\otimes_f$ defines a tensor structure on $\text{Mod}$. \end{thm}

\noindent The category $\text{Mod}$ together with the tensor product $\otimes_f$ is called quantum fusion tensor category (see \cite{ma} for the definition and complements on tensor categories). The aim of this section is to
prove this theorem. We warn that objects of the category
$\text{Mod}$ are not necessarily of finite length. Sometimes it is
required that the objects of a tensor category have a finite
composition series (for example see \cite{ce}), so $\text{Mod}$ is not a tensor category in this sense. However, we will prove in section \ref{flsection} that for a large class of quantum affinizations (including quantum affine
algebras and most quantum toroidal algebras) the subcategory of
finite length modules is stable under the monoidal bifunctor
$\otimes_f$, and so we get a tensor category in this sense.

{\it Proof of theorem \ref{bif}.} We already have proved the well-definedness in lemma \ref{proddef}. We have to show the associativity and existence of a neutral element. This will be formulated in the following two lemmas.

\begin{lem} $\otimes_f$ is associative.\end{lem}

\demo Let $r_1, r_2, r_3\geq 1$, $V_1\in\text{Mod}^{r_1}, V_2\in\text{Mod}^{r_2}, V_3\in\text{Mod}^{r_3}$. Let us prove that the identity defines an isomorphism between
the modules $V_1\otimes_f (V_2\otimes_f V_3)$ and $(V_1\otimes_f V_2)\otimes_f V_3$ as objects of $\text{Mod}^{r_1+r_2+r_3}$. The action of $g\in\U_q'(\hat{\Glie})$ is
given in the first case by : $(\text{Id}\otimes\Delta_{u^{r_2}})\circ\Delta_{u^{r_1}}$ and in the second case by :
$(\Delta_{u^{r_1}}\otimes\text{Id})\circ\Delta_{u^{r_1+r_2}}$. But it follows from lemma \ref{coass} that these maps are equal. If one $r_i$ is equal to $0$, $V_1\otimes_f (V_2\otimes_f V_3)\simeq (V_1\otimes_f V_2)\otimes_f V_3$ is a direct consequence of lemma \ref{extend}. 

\noindent The pentagon axiom is clearly satisfied as for the usual tensor category of vector spaces.\qed

\begin{lem}\label{neutral} $i(L_0)$ object of $\text{Mod}^0$ is a neutral object of $(\text{Mod},\otimes_f)$.\end{lem}

\demo It is a direct consequence of lemma \ref{extend} (the triangle axiom is clearly satisfied as for the usual tensor category of vector spaces).\qed

%\subsection{First applications}
%\subsubsection{Rationality of the fusion procedure}
%\noindent Note that some categories
%$\mathcal{O}^R(\U_q'(\hat{\Glie}))$ ($R\geq 0$) and a bilinear map
%$\otimes_R : \mathcal{O}(\U_q(\hat{\Glie}))\times
%\mathcal{O}^R(\U_q'(\hat{\Glie}))\rightarrow\mathcal{O}^{R+1}(\U_q'(\hat{\Glie}))$
%were considered in \cite{her04}. We have
%$\text{Mod}(\U_q(\hat{\Glie}))\subset\mathcal{O}^0(\U_q'(\hat{\Glie}))$,
%and $\otimes_0$ corresponds to the tensor product defined above for
%modules in $\text{Mod}(\U_q'(\hat{\Glie}))$, except that
%$V_1\otimes_0 V_2=(V_1\otimes_{\CC}V_2)\otimes \CC((u))$. So the
%lemma \ref{proddef} implies a rationality result :
%\begin{cor}\label{appun} Let $V_1,V_2\in\text{Mod}(\U_q(\hat{\Glie}))$. The subspace %$\CC(u)$-vector space $(V_1\otimes_{\CC} V_2)\otimes\CC(u)\subset (V_1\otimes_1 V_2)$ is %stable by the action of $\U_q'(\hat{\Glie})$.\end{cor}

\subsubsection{Application to $(l,u)$-highest weight simple modules of the category $\text{Mod}(\U_q'(\hat{\Glie}))$}
We use the notations of remark \ref{toref}. As an application of lemma \ref{proddef} we prove a "deformed version" of the "if" part of theorem \ref{intsimp}.

\begin{defi}\label{luweight} An $(l,u)$-weight $(\lambda,\Psi(u))$ is said to be dominant if for $i\in I$ there exists a polynomial $P_{i,u}(z)=(1-za_{i,1}u^{b_{i,1}})...(1-za_{i,N_i}u^{b_{i,N_i}})$ ($N_i\geq 0$, $a_{i,j}\in\CC^*$, $b_{i,j}\geq 0$) such that in $\CC[u^{\pm 1}][[z]]$ (resp. in $\CC[u^{\pm 1}][[z^{-1}]]$):
$${\sum}_{r\geq 0}\Psi_{i,\pm r}^{\pm}(u) z^{\pm r}= q_i^{N_i}\frac{P_{i,u}(zq_i^{-1})}{P_{i,u}(zq_i)}.$$\end{defi}

\noindent The notion of $(l,u)$-highest weight $\U_q'(\hat{\Glie})$-module is defined analogously to the notion of $l$-highest weight $\U_q(\hat{\Glie})$-module. By theorem \ref{afftrian} we get the existence of the simple module $L(\lambda,\Psi(u))$ of dominant $(l,u)$-highest weight $(\lambda,\Psi(u))$.

\begin{cor}\label{appdeux} For $(\lambda,\Psi(u))$ a dominant $(l,u)$-highest weight, we have $L(\lambda,\Psi(u))\in\text{Mod}(\U_q'(\hat{\Glie}))$.\end{cor}

\demo We see as in \cite[Lemma 32]{her04} that
$L(\lambda,\Psi(u))$ is the quotient of a tensor product in
$\text{Mod}$ of (fundamental) simple modules of
$\text{Mod}(\U_q(\hat{\Glie}))$. So the result follows from lemma \ref{proddef}.\qed

\noindent Note that this result in finer than \cite[Lemma 32]{her04} on the representation $L(\lambda,\Psi(u))$. The converse of this result is false (for example by using an evaluation morphism $\U_q'(\hat{sl_2})\rightarrow \U_q(sl_2)\otimes \CC(u)$, one can define a $2$-dimensional $(l,u)$-highest weight module with $(l,u)$-weight given by $P_u(z)=(1-z(u+u^{-1}))$ in the formula of definition \ref{luweight}).

\subsection{Example : fusion of Kirillov-Reshetikhin modules of type $sl_2$}\label{exkr}

The $\U_q(\hat{sl_2})$ Kirillov-Reshetikhin modules can be described explicitly because they can be realized as evaluation representations from simple finite dimensional representations of the quantum group of finite type $\U_q(sl_2)$. Let $r\geq 0$ and $a\in\CC^*$. The Kirillov-Reshetikhin module $W_r(a)$ is the $\U_q(\hat{sl_2})$-module of dimension $(r+1)$ with a basis $\{v_0,...,v_r\}$ and the action of the generators $x_r^+,x_r^-$ ($r\in\ZZ$), $\phi_{\pm m}^{\pm}$ ($m\geq 1$), $k^{\pm 1}$ given by :
$$x_r^+.v_j=(aq^{r-2j+2})^r[r-j+1]_qv_{j-1},$$
$$x_r^-.v_j=(aq^{r-2j})^r[j+1]_qv_{j+1},$$
$$k^{\pm 1}.v_j=q^{\pm (r-2j)}v_j,$$
$$\phi_{\pm m}^{\pm}.v_j=\pm (q-q^{-1})(q^{r-2j+1}a)^{\pm m}(q^{\mp m}[r-j]_q[j+1]_q-q^{\pm m}[j]_q[r-j+1]_q)v_j,$$
where $0\leq j\leq r$ and we denote $v_{-1}=v_{r+1}=0$.
In particular we have as a power series in $z$ (resp. in $z^{-1}$) :
$$\phi^{\pm}(z).v_j=q^{r-2j}\frac{(1-zaq^{-r})(1-zaq^{r+2})}{(1-zaq^{r-2j+2})(1-zaq^{r-2j})}v_j.$$

\noindent Consider $W_r(a), W_{r'}(b)$ two Kirillov-Reshetikhin
modules with respective basis $\{v_0,...,v_r\}$,
\\$\{v_0',...,v_{r'}'\}$. One can describe explicitly the fusion
module $i(W_r(a))\otimes_f i(W_{r'}(b))$. It has a $\CC(u)$-basis
$\{v_j\otimes v_k'|0\leq j\leq r, 0\leq k\leq
r'\}$, and the action of $\U_q'(\hat{sl_2})$ is given by
\begin{equation}\label{krplus}x_m^+.(v_j\otimes v_k')=(aq^{r-2j+2})^m
\alpha_{j,k}(v_{j-1}\otimes v_k') +(ubq^{r'-2k+2})^m \beta_{j,k}
(v_j\otimes v_{k-1}'),\end{equation}
\begin{equation}\label{krxmoins}x_m^-.(v_j\otimes v_k')=(ubq^{r'-2k})^m
\gamma_{j,k}(v_{j}\otimes v_{k+1}')
+(aq^{r-2j})^m\mu_{j,k}(v_{j+1}\otimes v_{k}'),\end{equation}
\begin{equation}\label{eigenphi}\phi^{\pm}(z).(v_j\otimes
v_k')=q^{r-2jr'-2k}\frac{(1-zaq^{-r})(1-zaq^{r+2})(1-zubq^{-r'})(1-zubq^{r'+2})}{(1-zaq^{r-2j+2})(1-zaq^{r-2j})(1-zubq^{r'-2k+2})(1-zubq^{r-2k})}(v_j\otimes
v_k'),\end{equation} 
where $\alpha_{j,k}=[r-j+1]_q$,
$\gamma_{j,k}=[k+1]_q$, and $\beta_{j,k}$ is equal to 
$$[r'-k+1]_q(q^{2j-r}+(q^{-1}-q)ua^{-1}bq^{r'-2k+1-r+2j}(\frac{q[r-j]_q[j+1]_q}{1-ua^{-1}bq^{r'-2k+2-r+2j}}-\frac{q^{-1}[j]_q[r-j+1]_q}{1-ua^{-1}bq^{r'-2k-r+2j}})),$$
and $\mu_{j,k}$ is equal to  $$[j+1]_q(q^{r'-2k}+(q-q^{-1})ua^{-1}bq^{r'-2k+1-r+2j}(\frac{q^{-1}[r'-k]_q[k+1]_q}{1-ua^{-1}bq^{r'-2k-r+2j}}-\frac{q
[k]_q[r-k+1]_q}{1-ua^{-1}bq^{r'-2k-r+2j+2}})).$$ 
Note that $\alpha_{j,k}, \beta_{j,k}, \gamma_{j,k}, \mu_{j,k}$ are independent of $m$.

\noindent (Observe that on this example the action is
rational, as proved in lemma \ref{proddef}). 

\begin{rem}\label{usucop} For the usual tensor product of quantum affine algebras, certain tensor products of $l$-highest weight modules are $l$-highest weight (see \cite[Theorem 4]{c} and \cite[Theorem 9.1]{kas}) but not all. For example for $\U_q(\hat{sl_2})$, $L_a\otimes L_{aq^2}$ is $l$-highest weight but $L_{aq^2}\otimes L_a$ is not $l$-highest weight (where the $L_b$ are the fundamental representations defined in section \ref{int}).\end{rem}

\noindent In contrast to this well-known situation, in the following examples the fusion modules are always of $l$-highest weight (we will see in theorem \ref{lhw} that this observation is a particular case of a more general picture for the quantum fusion tensor category).

\begin{prop}\label{luhw} The fusion module $V=i(W_r(a))\otimes_f i(W_{r'}(b))$ is a simple $(l,u)$-highest weight module.\end{prop}

\demo First let us prove that $V$ is of $(l,u)$-highest weight. Let $W$ be the sub $\U_q'(\hat{\Glie})$-module of $V$ generated by $v_0\otimes v_0'$. It suffices to prove by induction on $K\geq 0$ that $((j+k=K)\Rightarrow (v_j\otimes v_k'\in W))$. For $K=0$ it is clear. Let $K\geq 1$ and $j,k\geq 0$ satisfying $j+k=K-1$. By definition, $\gamma_{j,k}, \mu_{j,k}$ are not equal to zero. As $ubq^{r'-2k}\neq aq^{r-2j}$, equation (\ref{krxmoins}) implies that $v_{j+1}\otimes v_k'$ and $v_j\otimes v_{k+1}'$ are in ${\sum}_{r\in\ZZ}\CC(u)x_r^-.(v_j\otimes v_k')\subset W$. So $\sum_{\{(a,b)|a+b=K\}}\CC(v_a\otimes v_b')\subset W$.

\noindent Let us prove that $V$ is simple. Suppose that $V'$ is a
proper submodule of $V$. Suppose that the eigenvalues given in equation (\ref{eigenphi})
for $v_{j_1}\otimes v_{k_1}'$ and $v_{j_2}\otimes v_{k_2}'$  are equal. Then we have $$(1-zaq^{r-2j_1+2})(1-zaq^{r-2j_1})=(1-zaq^{r-2j_2+2})(1-zaq^{r-2j_2})$$ 
and 
$$(1-zubq^{r'-2k_1+2})(1-zubq^{r-2k_1})=(1-zubq^{r'-2k_2+2})(1-zubq^{r-2k_2}).$$ 
So $j_1=j_2$
and $k_1=k_2$ as the conditions $aq^{r-2j_1}=aq^{r-2j_2+2}$ and $aq^{r-2j_1+2}=aq^{r-2j_2}$ (resp. $aq^{r-2k_1}=aq^{r-2k_2+2}$ and $aq^{r-2k_1+2}=aq^{r-2k_2}$) can not be simultaneously satisfied. So the operators $\phi^{\pm}_{\pm m}$ have a diagonalizable action on $V$ with common $l$-weight spaces of dimension $1$. As $V'$ is stable under the action of the $\phi^{\pm}_{\pm m}$, it is of the form $V'={\bigoplus}_{(j,k)\in J}\CC(u) (v_j\otimes v_k')$ where $J\subset \{0,...,r\}\times\{0,...,r'\}$. Consider $v_j\otimes v_k'\in V'$ such that $j+k$ is minimal for this property. If $V'\neq V$, we have $j+k > 0$ and $V'\cap (\bigoplus_{\{(s,t)|s+t < j+k\}}\CC(u)(v_{s}\otimes v_{t}'))=\{0\}$. But from equation (\ref{krplus}), we have $x^+(z).(v_j\otimes v_k')\neq 0$, contradiction. \qed

\noindent The results of this section will be used in more detail in section \ref{tsyst}.

\section{$\mathcal{A}$-forms and specializations}\label{saform}

The aim of this section is to explain how to go back from $\text{Mod}$ to the usual category $\text{Mod}(\U_q(\hat{\Glie}))$. Consider the ring $\mathcal{A}=\{\frac{f(u)}{g(u)}|f(u),g(u)\in\CC[u^{\pm}], g(1)\neq 0\}\subset\CC(u)$. We prove the existence of certain $\mathcal{A}$-forms (definition \ref{aform}) and we study their specialization at $u=1$. The idea of considering $\mathcal{A}$-form is inspired from crystal basis theory (see \cite{kas0, lu0}), but instead of using $q$ as a deformation parameter in the ring $\mathcal{A}$, we use the deformation parameter $u$ of the Drinfeld coproduct. We prove that one can define $\mathcal{A}$-forms of cyclic modules of the category $\text{Mod}(\U_q'(\hat{\Glie}))$ (theorem \ref{aformcy}). It gives rise to the following two constructions : one can construct a fusion $l$-highest weight module from two $l$-highest weight modules (definition \ref{defsf}), and in the case of quantum affine algebras one can define a fusion bifunctor for the category of finite dimensional representations (theorem \ref{defsd}).

\subsection{$\mathcal{A}$-forms} In this section we recall standard definitions and properties of $\mathcal{A}$-forms which are some $\mathcal{A}$-"lattice" of $\CC(u)$-vector spaces. We also study explicit examples.

\noindent $\mathcal{A}$ is a local
principal ring. Let
$\U_q^u(\hat{\Glie})=\U_q(\hat{\Glie})\otimes_{\CC}\mathcal{A}\subset\U_q'(\hat{\Glie})$.
The algebras $\U_q^{u,+}(\hat{\Glie}), \U_q^{u,-}(\hat{\Glie}),
\U_q^u(\hat{\Hlie})\subset\U_q^u(\hat{\Glie})$ are defined in an
obvious way.

\subsubsection{Definition and first properties} 

\begin{defi}\label{aform} Let $V\in\text{Mod}(\U_q'(\hat{\Glie}))$. An $\mathcal{A}$-form $\tilde{V}$ of $V$ is a sub $\mathcal{A}$-module of $V$ satisfying :

i) $\tilde{V}$ is stable under $\U_q^u(\hat{\Glie})$,

ii) $\tilde{V}$ generates the $\CC(u)$-vector space $V$,

iii) for $\lambda\in \Hlie^*$, $\tilde{V}\cap V_{\lambda}$ is a finitely generated $\mathcal{A}$-module.\end{defi}

\noindent Note that as $\mathcal{A}$ is principal and $V$ is torsion free, the property iii) implies that $\tilde{V}\cap V_{\lambda}$ is free. In particular it makes sense to consider the rank and basis of $\tilde{V}\cap V_{\lambda}$.

\noindent This notion of $\mathcal{A}$-form is analog to the lattice $L$ considered in crystal theory \cite{kas0, lu0}. The properties discussed in this subsection are standard (see \cite{lu, kas3} for general results on $\mathcal{A}$-forms). We give proofs for the convenience of the reader.

\begin{lem}\label{firstprop} Let $V\in\text{Mod}(\U_q'(\hat{\Glie}))$ and $\tilde{V}$ be an $\mathcal{A}$-form of $V$. Then 

1) $V\simeq \tilde{V}\otimes_{\mathcal{A}}\CC(u)$ as $\U_q'(\hat{\Glie})$-modules,

2) for all $\lambda\in\Hlie^*$, we have $\text{rk}_{\mathcal{A}}(\tilde{V}\cap V_{\lambda})=\text{dim}_{\CC(u)}(V_{\lambda})$,

3) for all $v\in V$, there is a unique $n(v)\in\ZZ$ such that $(1-u)^{n(v)}v\in \tilde{V}$ and $(1-u)^{n(v)-1}v\notin\tilde{V}$,

4) let $W\subset V$ be a $\U_q^u(\hat{\Glie})$-submodule of $V$. Then $W\cap\tilde{V}$ is an $\mathcal{A}$-form of $W$.\end{lem}

\demo 1) and 2) : for all $\lambda\in \Hlie^*$, a basis of $\tilde{V}\cap V_{\lambda}$ as an $\mathcal{A}$-module is also a basis of the $\CC(u)$-vector space
$V_{\lambda}$.

3) It suffices to consider the case $v\in V_{\lambda}$. Let us write it $v_{\lambda}={\sum}_{k=1,\cdots,K}f_i(u) w_i$ where $\{w_i\}_{i=1,\cdots,K}$ is an
$\mathcal{A}$-basis of $V_{\lambda}\cap\tilde{V}$. Then $n(v)$ is the unique $n(v)\in\ZZ$ satisfying : $\forall i\in \{1,...,K\}$, $(1-u)^{n(v)}f_i(u)\in\mathcal{A}$
and there exists $i\in \{1,...,K\}$ such that $(1-u)^{n(v)-1}f_i(u)\notin\mathcal{A}$.

4) Denote $\tilde{W}=W\cap\tilde{V}$, and let us check the properties of definition \ref{aform} : i) is clear. For ii), let $v\in W$. There is $P\in \CC[u^{\pm}]-\{0\}$ such that $v'=P v\in\tilde{V}$. So $v'\in\tilde{W}$ and $\tilde{W}$ generates the $\CC(u)$-vector space $W$. For the property $iii)$, we have $\tilde{W}\cap W_{\lambda}=V_{\lambda}\cap W\cap\tilde{V}\subset \tilde{V}\cap V_{\lambda}$ is a finitely generated $\mathcal{A}$-module.\qed

\noindent We have directly from definition \ref{aform} :

\begin{lem}\label{sumaform} Let $V\in\text{Mod}(\U_q'(\hat{\Glie}))$ and $\tilde{V}\subset V$ (resp. $\tilde{V'}\subset V$) be an $\mathcal{A}$-form of $\U_q'(\hat{\Glie}).\tilde{V}$ (resp. of $\U_q'(\hat{\Glie}).\tilde{V'}$). Then $\tilde{V}+\tilde{V'}$ is an $\mathcal{A}$-form of $\U_q'(\hat{\Glie}).(\tilde{V}+\tilde{V'})$.\end{lem}

\subsubsection{Specialization} One can consider the specialization of an $\mathcal{A}$-form at $u=1$ :

\begin{defi} Let $V\in\text{Mod}(\U_q'(\hat{\Glie}))$ and $\tilde{V}$ an $\mathcal{A}$-form of $V$. Then we denote by $(\tilde{V})_{u=1}$ the $\U_q(\hat{\Glie})$-module $\tilde{V}/((u-1)\tilde{V})$.\end{defi}

\begin{lem} We have $(\tilde{V})_{u=1}\in\text{Mod}(\U_q(\hat{\Glie}))$.\end{lem}

\demo It follows from lemma \ref{firstprop} that  $$\text{dim}_{\CC}(((\tilde{V})_{u=1})_{\lambda})=\text{rk}_{\mathcal{A}}(\tilde{V}\cap V_{\lambda})=\text{dim}_{\CC(u)}(V_{\lambda}).$$
So $(\tilde{V})_{u=1}\in\text{Mod}(\U_q(\hat{\Glie}))$.\qed

\noindent In general the specialization at $u=1$ of two $\mathcal{A}$-forms of an $\U_q'(\hat{\Glie})$-module are not necessarily isomorphic, as illustrated in the following examples.

\subsubsection{Examples}\label{exspe} Let $L_1=L_{1,1}=\CC v_0\oplus\CC v_1$ and $L_2=L_{1,q^2}=\CC w_0\oplus \CC w_1$ be two
$\U_q(\hat{sl_2})$-fundamental representations. The fusion modules
$V=i(L_1)\otimes_f i(L_2)$ and $V'=i(L_2)\otimes_f i(L_1)$ can be
 described explicitly (see section \ref{exkr}; they have also been
studied in \cite[section 6.6]{her04} with a different formalism). Consider 
the $\CC(u)$-basis of $V$ (resp. of $V'$): $f_0=v_0\otimes
w_0, f_1=v_1\otimes w_0, f_2=v_0\otimes w_1, f_3=v_1\otimes w_1$
(resp. $f_0'=w_0\otimes v_0, f_1'=w_1\otimes v_0, f_2'=w_0\otimes
v_1, f_3'=w_1\otimes v_1$). The action of $\U_q'(\hat{sl_2})$ on $V$
is given by :
$$
\begin{array}{l|l|l|l|l}
   & f_0& f_1& f_2 & f_3\\ \hline
 x_r^+ & 0 & f_0& q^{2r-1}u^r \frac{1-q^4u}{1-uq^2} f_0& f_2+q^{1+2r}u^r\frac{1-u}{1-uq^2}f_1\\ \hline
x_{r}^- & u^rq^{2r}f_2+q\frac{1-u}{1-q^2u}f_1& u^rq^{2r}f_3& q^{-1}\frac{1-q^4u}{1-q^2u}f_3& 0 \\ \hline
\phi^{\pm}(z) & q^2\frac{(1-q^{-2}z)(1-uz)}{(1-z)(1-q^2uz)}f_0& \frac{(1-q^2z)(1-uz)}{(1-z)(1-q^2 uz)}f_1& \frac{(1-q^{-2}z)(1-q^4uz)}{(1-z)(1-q^2uz)}f_2& q^{-2}\frac{(1-q^2 z)(1-q^4uz)}{(1-z)(1-q^2 uz)}f_3\end{array}.$$
The action of $\U_q'(\hat{sl_2})$ on $V'$ is given by :
$$
\begin{array}{l|l|l|l|l}
   & f'_0& f'_1& f'_2 & f'_3\\ \hline
 x_r^+ & 0 & q^{2r}f'_0& q^{-1}u^r \frac{1-u}{1-uq^{-2}} f'_0& q^{2r}f'_2+q u^r\frac{1-q^{-4}u}{1-uq^{-2}}f'_1\\ \hline
x_{r}^- & u^rf'_2+q^{2r+1}\frac{1-q^{-4}u}{1-q^{-2}u}f'_1& u^r f'_3& q^{2r-1}\frac{1-u}{1-q^{-2}u}f'_3& 0 \\ \hline
\phi^{\pm}(z) & q^2\frac{(1-z)(1-q^{-2}uz)}{(1-q^2z)(1-uz)}f'_0& \frac{(1-q^4z)(1-q^{-2}uz)}{(1-q^2z)(1-uz)}f'_1& \frac{(1-z)(1-q^2uz)}{(1-q^2z)(1-uz)}f'_2& q^{-2}\frac{(1-q^4 z)(1-q^2uz)}{(1-q^2z)(1-uz)}f'_3\end{array}.$$
It follows from proposition \ref{luhw} that $V$ and $V'$ are simple $(l,u)$-highest weight modules.

\noindent In particular
$$\tilde{V}=\U_q^u(\hat{\Glie}).f_0=\mathcal{A}.f_0\oplus \mathcal{A} (1-u).f_1\oplus \mathcal{A}.f_2\oplus \mathcal{A} f_3,$$
$$\tilde{V'}=\U_q^u(\hat{\Glie}).f_0'=\mathcal{A}.f'_0\oplus \mathcal{A} f'_1\oplus \mathcal{A}.f'_2\oplus \mathcal{A} f'_3,$$
are respectively $\mathcal{A}$-forms of $V$ and $V'$. We use the same notations for the vectors in $\tilde{V}$ and in $(\tilde{V})_{u=1}$ (resp. in $\tilde{V}'$ and in $(\tilde{V}')_{u=1}$). The $\U_q(\hat{sl_2})$-modules
$(\tilde{V})_{u=1}$ and $(\tilde{V'})_{u=1}$ are $l$-highest weight
modules, but not simple :  $\CC . (1-u)f_1$ (resp. $\CC.f_2'$) is a $\U_q(\hat{sl_2})$-submodule of $(\tilde{V})_{u=1}$
(resp. of $(\tilde{V'})_{u=1}$) of dimension $1$.

\noindent Note that we have an isomorphism of
$\U_q(\hat{sl_2})$-modules $\sigma : (\tilde{V})_{u=1}\rightarrow
(\tilde{V'})_{u=1}$ defined by $\sigma(f_0)=f_0'$,
$\sigma((1-u)f_1)=(q^{-1}-q)f_2'$,
$\sigma(f_2)=(q+q^{-1})f_1'$, $\sigma(f_3)=f_3'$.

\noindent Consider the following respective $\mathcal{A}$-forms of $V$ and $V'$ :
$$\tilde{V}_f={\sum}_{i=0,1,2,3}\U_q^u(\hat{sl_2}).f_i=\mathcal{A}f_0\oplus\mathcal{A}f_1\oplus\mathcal{A}f_2\oplus\mathcal{A}f_3\supsetneq \tilde{V},$$
$${\sum}_{i=0,1,2,3}\U_q^u(\hat{sl_2}).f_i'=\mathcal{A}f_0'\oplus\mathcal{A}f_1'\oplus\mathcal{A}f_2'\oplus\mathcal{A}f_3'=\tilde{V'}.$$
We use the same notation for the vectors in $\tilde{V}_f$ and in $(\tilde{V}_f)_{u=1}$. $(\tilde{V}_f)_{u=1}$ is not an $l$-highest weight module, is cyclic generated by $f_1$ and has a submodule of dimension $3$, namely $\CC f_0\oplus \CC f_2\oplus \CC f_3$. So $(\tilde{V}_f)_{u=1}$ and $(\tilde{V})_{u=1}$ are not isomorphic.

\subsection{$\mathcal{A}$-form of cyclic modules} In this section we study $\mathcal{A}$-forms of cyclic modules which will be used later (in particular we will study the interesting properties of their specializations in other sections). The main result of this section is :

\begin{thm}\label{aformcy} Let $V\in\text{Mod}(\U_q'(\hat{\Glie}))$ and $v\in V-\{0 \}$. Then $\tilde{V}(v)=\U_q^u(\hat{\Glie}).v$ is an $\mathcal{A}$-form of $\U_q'(\hat{\Glie}).v$.
Moreover $(\tilde{V}(v))_{u=1}\in\text{Mod}(\U_q(\hat{\Glie}))$ is a non zero cyclic $\U_q(\hat{\Glie})$-module generated by $v$.\end{thm}

\noindent This theorem is proved in this section. We can suppose that $V$ is a non zero cyclic $\U_q'(\hat{\Glie})$-module generated by $v$.

\begin{lem}\label{hfree} Let $V\in\text{Mod}(\U_q'(\hat{\Glie}))$ and $F\subset V$ be a finitely generated $\mathcal{A}$-submodule of $V$. Then $\mathcal{A}.\U_q(\hat{\Hlie}).F$ is a finitely generated $\mathcal{A}$-module.\end{lem}

\demo Let $W\subset V$ as in definition \ref{ucat} and write
$F={\sum}_{j=1,...,m}\mathcal{A}.f_j$ where $f_j\in V$. For
each $j\in\{1,...,m\}$, we can write a finite sum
$f_j={\sum}_k\Psi_{j,k}f_{j,k}$ where $f_{j,k}\in
W_{\lambda_{j,k}}$ ($\lambda_{j,k}\in\Hlie^*$) and
$\Psi_{j,k}\in\CC(u)$. But
$\mathcal{A}.\U_q(\hat{\Hlie})f_{j,k}\subset
\mathcal{A}.W_{\lambda_{j,k}}$, and so
$\U_q(\hat{\Hlie}).F\subset{\sum}_{j,k}\Psi_{j,k}\mathcal{A}.W_{\lambda_{j,k}}$
is a finitely generated $\mathcal{A}$-module.\qed

\noindent Let us study the case of $(l,u)$-highest weight modules, which is a first step in the proof of theorem \ref{aformcy} :

\begin{lem}\label{aformhw} Let $V$ be an $(l,u)$-highest weight module in the category $\text{Mod}(\U_q(\hat{\Glie}))$ and let $v$ be an highest weight vector. Then
$\tilde{V}=\U_q^u(\hat{\Glie}).v$ is an $\mathcal{A}$-form of $V$.\end{lem}

\demo Properties i), ii) of definition \ref{aform} are clear. Let us prove property iii) : let $\lambda$ be the weight of $v$. For $\mu$ a weight of $V$, we
have $\mu\in\lambda-Q^+$. Let us prove the result by induction on $h'(\mu)=h(\lambda-\mu)$. For $h'(\mu)=0$ it is clear, and in general let us prove that $\Lambda_{l+1}={\sum}_{\{\mu|h'(\mu)=l+1\}}\tilde{V}\cap V_{\mu}$ is a finitely generated $\mathcal{A}$-module. But we have $\Lambda_{l+1}={\sum}_{\{i\in I,
m\in\ZZ, \mu|h'(\mu)=l\}}x_{i,m}^-.(V_{\mu}\cap\tilde{V})$. It follows from formula (\ref{actcartld}) that for $i\in I, m\neq 0$, we have
$x_{i,m}^-=\frac{-m}{[2r_i]_q}(h_{i,m}x_{i,0}^- -x_{i,0}^-h_{i,m})$. Moreover we have $\U_q(\hat{\Hlie}).(V_{\mu}\cap\tilde{V})=V_{\mu}\cap\tilde{V}$. In particular
$\Lambda_{l+1}\subset{\sum}_{\{i\in I,\mu|h'(\mu)=l\}}\U_q(\hat{\Hlie}).x_{i,0}^-(V_{\mu}\cap \tilde{V})$. We can conclude with lemma \ref{hfree}.\qed

\noindent Note that the rationality in $u$ was a crucial point of this proof.

\noindent (Note that in \cite{her04} we considered
$\CC[u^{\pm}]\U_q(\hat{\Glie}).v$ for a simple $(l,u)$-highest weight
module $L$ and called it a $\CC[u^{\pm}]$-form of $L$. As
$\U_q^u(\hat{\Glie}).v=\mathcal{A}.\CC[u^{\pm}]\U_q(\hat{\Glie}).v$,
it is a particular case of the point of view of this paper.)

{\it Proof of theorem \ref{aformcy} :}

\noindent Let us prove that we get an $\mathcal{A}$-form : properties i), ii) of definition \ref{aform} are clear. Let us check property iii). First consider
$$d(v)=\text{Max}\{h(\text{wt}(v')-\text{wt}(v))|v'\in V\text{ and }\text{wt}(v')-\text{wt}(v)\in Q^+\}.$$
Let us prove by induction on $d(v)\geq 0$,
that for all $\lambda\in\Hlie^*$, $\U_q^u(\hat{\Glie}).v\cap (\U_q'(\hat{\Glie}).v)_{\lambda}$ is a finitely generated $\mathcal{A}$-module. For $d(v)=0$, the result is
proved as in lemma \ref{aformhw}. In general : the triangular decomposition of $\U_q'(\hat{\Glie})$ (theorem \ref{afftrian}) gives $\U_q^u(\hat{\Glie}).v=A+B$ where
$A=\U_q^{u,-}(\hat{\Glie})\U_q^u(\hat{\Hlie}).v$ and $B={\sum}_{i\in I,m\in\ZZ}\U_q^u(\hat{\Glie}).x_{i,m}^+.v$. Moreover for $\lambda\in\Hlie^*$,
$(\U_q^u(\hat{\Glie}).v)\cap V_{\lambda}=A\cap V_{\lambda}+B\cap V_{\lambda}$. We see as in the proof of lemma \ref{aformhw} that $A\cap V_{\lambda}$ is a finitely
generated $\mathcal{A}$-module. It follows from formula (\ref{actcartld}) that for $i\in I, m\neq 0$, we have $x_{i,m}^+=\frac{m}{[2r_i]_q}(h_{i,m}x_{i,0}^+ -x_{i,0}^+h_{i,m})$. In
particular $B\subset {\sum}_{i\in I}\U_q^u(\hat{\Glie}).x_{i,0}^+.\U_q^u(\hat{\Hlie}).v$. It follows from lemma \ref{hfree} that $\U_q^u(\hat{\Hlie}).v$ is
a finitely generated $\mathcal{A}$-module. For $i\in I$, and $v'\in x_{i,0}^+.\U_q^u(\hat{\Hlie}).v$, we have $d(v')<d(v)$, and so it follows from the induction hypothesis that
for all $\lambda\in\Hlie^*$, $\U_q^u(\hat{\Glie}).v'\cap V_{\lambda}=\U_q^u(\hat{\Glie}).v'\cap (\U_q^u(\hat{\Glie}).v')_{\lambda}$ is a finitely generated
$\mathcal{A}$-module, and so $B\cap V_{\lambda}$ is a finitely generated
$\mathcal{A}$-module.

\noindent Let us prove that the image of $v$ in $(\tilde{V}(v))_{u=1}$ generates
$(\tilde{V}(v))_{u=1}$ : first we have $v\notin (1-u)\tilde{V}(v)$.
Indeed if there is an $r<0$ such that $(1-u)^rv\in \tilde{V}(v)$, then
for all $l\geq 1$, $(1-u)^{lr}v\in\tilde{V}(v)$ and so if
$v\in{\sum}_{\lambda\in\Lambda}\tilde{V}(v)\cap V_{\lambda}$
($\Lambda$ finite), there is $\lambda\in\Lambda$ such that
$\tilde{V}(v)\cap V_{\lambda}$ is not a finitely generated
$\mathcal{A}$-module, contradiction. In particular
$(\tilde{V}(v))_{u=1}\neq\{0\}$. Next let
$w\in \tilde{V}(v)-((1-u)\tilde{V}(v))$ and
$G(u)\in\U_q^u(\hat{\Glie})$ such that $w=G(u).v$. Then in
$(\tilde{V}(v))_{u=1}$, we have $w=G(1).v$ (where $G(1)$ is the image of $G(u)$ by the projection
$\U_q^u(\hat{\Glie})\rightarrow\U_q(\hat{\Glie})$).\qed

\subsection{Fusion of $l$-highest weight $\U_q(\hat{\Glie})$-modules} In this section we explain how to construct fusion of $l$-highest weight $\U_q(\hat{\Glie})$-modules by using the quantum fusion tensor category. We use the notations of remark \ref{toref}. 

\noindent We get the following first application of theorem \ref{aformcy} :

\begin{cor}\label{infosf} Let $r\geq 1$ and $V_1, \cdots,V_r$ be $l$-highest weight modules in the category $\text{Mod}(\U_q(\hat{\Glie}))$. Consider the fusion module
$W=i(V_1)\otimes_f i(V_2)\otimes_f \cdots\otimes_f i(V_r)\in\text{Mod}(\U_q'(\hat{\Glie}))$. Let $v_1\in V_1,\cdots,v_r\in V_r$ be highest weight vectors and $(P_i^1)_{i\in
I},...,(P_i^r)_{i\in I}$ the respective Drinfeld polynomials. Then

1) The $\U_q(\hat{\Glie})$-module $(\tilde{W}(v))_{u=1}\in\text{Mod}(\U_q(\hat{\Glie}))$ is an $l$-highest weight module of $l$-highest weight vector $v=v_1\otimes \cdots\otimes v_r$. 

2) The Drinfeld polynomials of the $l$-highest weight of $(\tilde{W}(v))_{u=1}$ are $(P_i^1(z)P_i^2(z)\cdots P_i^r(z))_{i\in I}$.\end{cor}

\demo 1) It follows from theorem \ref{aformcy} that $(\tilde{W}(v))_{u=1}$ is cyclic generated by $v$. Moreover $v_1,\cdots,v_r$ are $l$-highest weight vectors, and so formulas (\ref{deltaup}), (\ref{deltauh}), (\ref{deltauh2}) give that $v$ is an $l$-highest weight vector. 

2) Clear from equations (\ref{deltauh}) and (\ref{deltauh2}).\qed

\begin{defi}\label{defsf} The module $(\tilde{W}(v))_{u=1}$ is denoted by $V_1*_f V_2*_f\cdots *_f V_r$ and is called the fusion module of $V_1,V_2,\cdots,V_r$.\end{defi}

\noindent Examples : in section \ref{exspe}, for the $\U_q(\hat{sl_2})$-$l$-highest weight modules $L_1$, $L_2$ we computed explicitly $L_1*_f L_2=(\tilde{V})_{u=1}$ and $L_2*_f L_1=(\tilde{V}')_{u=1}$ and we defined an isomorphism $\sigma : L_1*_f L_2\simeq L_2*_f L_1$. Moreover we noticed that $L_1*_f L_2$ is not semi-simple.

\noindent Other examples and applications will be studied in the section \ref{tsyst}.

\subsection{Fusion of finite dimensional representations}\label{tensd} We use the notations of remark \ref{toref}. In this subsection we study another application of the quantum fusion tensor category : it allows to
define a bifunctor on the category of finite dimensional representations of $\U_q(\hat{\Glie})$. This bifunctor will not be used in the rest of this paper, we hope to
study it in more details in another paper.

\noindent We recall that $i$ is a functor from $\text{Mod}(\U_q(\hat{\Glie}))$ to $\text{Mod}(\U_q'(\hat{\Glie}))$ (see proposition \ref{inclu}).

\begin{cor}\label{aideod} Let $V_1, V_2$ be two finite dimensional representations of $\U_q(\hat{\Glie})$. Then $$\U_q^u(\hat{\Glie}).(V_1\otimes_{\CC} V_2)\subset i(V_1)\otimes i(V_2)$$
is an $\mathcal{A}$-form of $i(V_1)\otimes i(V_2)$.\end{cor}

\demo Let $\{v_{\alpha}\}_{1\leq \alpha\leq p}, \{w_{\beta}\}_{1\leq\beta\leq p'}$ be $\CC$-basis respectively of $V_1$ and $V_2$. We have
\begin{equation}\label{fsum}\U_q^u(\hat{\Glie}).(V_1\otimes_{\CC} V_2)={\sum}_{1\leq \alpha\leq p , 1\leq\beta\leq p'}(\U_q^u(\hat{\Glie}).(v_{\alpha}\otimes w_{\beta}))\subset (V_1\otimes_{\CC}V_2)\otimes_{\CC} \CC(u).\end{equation} 
From theorem \ref{aformcy}, each $\U_q^u(\hat{\Glie}).(v_{\alpha}\otimes w_{\beta})$ is an $\mathcal{A}$-form of $\U_q'(\hat{\Glie}).(v_{\alpha}\otimes w_{\beta})$. As the sum in equation (\ref{fsum}) is finite, we can conclude with lemma \ref{sumaform}.\qed

Because of corollary \ref{aideod}, it makes sense to define :

\begin{defi} We denote by $V_1 \otimes_d V_2$ the $\U_q(\hat{\Glie})$-module $(\U_q^u(\hat{\Glie}).(V_1\otimes_{\CC} V_2))_{u=1}$.\end{defi}

\noindent {\it Examples} : in section \ref{exspe}, for the $\U_q(\hat{sl_2})$ $l$-highest weight modules $L_1$ and $L_2$ we computed explicitly $L_1\otimes_d
L_2=(\tilde{V}_f)_{u=1}$ and $L_2\otimes_d L_1=(\tilde{V}_f')_{u=1}$. Note that $L_1\otimes_d L_2$ is not isomorphic to $L_2\otimes_d L_1$, and that $L_1\otimes_d L_2$
is not isomorphic to $L_1*_f L_2$. Note that in general if $V_1$ and $V_2$ are semi-simple, $V_1\otimes_d V_2$ is not necessarily semi-simple.

\noindent Let $\text{Modf}(\U_q(\hat{\Glie}))$ be the subcategory of finite dimensional representations in $\text{Mod}(\U_q(\hat{\Glie}))$. If $\U_q(\hat{\Glie})$ is a
quantum affine algebra, the simple integrable $l$-highest weight modules are objects of $\text{Mod}(\U_q(\hat{\Glie}))$ (see \cite{Cha2}).

\begin{thm}\label{defsd} $\otimes_d$ defines a bifunctor $\otimes_d : \text{Modf}(\U_q(\hat{\Glie}))\times \text{Modf}(\U_q(\hat{\Glie}))\rightarrow \text{Modf}(\U_q(\hat{\Glie}))$.\end{thm}

\demo As $i(V_1)\otimes_f i(V_2)$ is a finite dimensional $\CC(u)$-vector space, $V_1\otimes_d V_2$ is a finite dimensional $\CC$-vector space, and so necessarily is an object of $\text{Modf}(\U_q(\hat{\Glie}))$. Consider $V_1, V_2, V_1', V_2'$ objects of $\text{Modf}(\U_q(\hat{\Glie}))$ and $f_1:V_1\rightarrow V_1'$, $f_2:V_2\rightarrow V_2'$ two morphisms of $\U_q(\hat{\Glie})$-module. From theorem \ref{bif} we have a morphism of $\U_q'(\hat{\Glie})$-module $f_1\otimes_f f_2: i(V_1)\otimes_f i(V_2)\rightarrow i(V_1')\otimes_f i(V_2')$. As $(f_1\otimes_f f_2)(\U_q^u.(V_1\otimes_{\CC}V_2))\subset \U_q^u.(V_1'\otimes_{\CC}V_2')$ and $(f_1\otimes_f f_2)((1-u)\U_q^u.(V_1\otimes_{\CC}V_2))\subset (1-u)\U_q^u.(V_1'\otimes_{\CC}V_2')$, we get a morphism $f_1\otimes_d f_2 : V_1\otimes_d V_2\rightarrow V_1'\otimes_d V_2'$.\qed

\noindent For quantum affine algebras, it should be interesting to
relate $\otimes_d$ to the usual tensor category structure.

\section{Finite length property and quantum fusion tensor category}\label{flsection} In this section we prove that for a large class of quantum affinizations (including
quantum affine algebras and most quantum toroidal algebras) the subcategory of finite length modules is stable under the monoidal bifunctor $\otimes_f$ (theorem
\ref{prectens}). In particular we get a tensor category in the sense of \cite{ce} (where the finite length property is required). We construct and use for this purpose a
generalization of the Frenkel-Reshetikhin morphism of $q$-characters \cite{Fre}
for the category $\text{Mod}$ and prove that it is compatible with
the monoidal structure (these results are valid for all quantum
affinizations). In particular it gives a representation theoretical
interpretation of the combinatorial fusion product defined in
\cite{her04} for the semi-simplified category (corollary \ref{inter}).

\noindent First let us recall the following consequence of theorem \ref{intsimp} :

\begin{cor}\label{simpmod} The simple modules of the category $\text{Mod}(\U_q(\hat{\Glie}))$ are the simple $l$-highest weight modules $L(\lambda,\Psi)$ where $(\lambda,\Psi)$ is a dominant $l$-weight.\end{cor}

\demo As a module in $\text{Mod}(\U_q(\hat{\Glie}))$ has at least one $l$-highest weight vector, the result is a direct consequence of theorem \ref{intsimp}.\qed

\subsection{Generalized $q$-characters and fusion product \cite{her04}} The theory of $q$-characters was introduced in \cite{Fre} for quantum affine algebras. In this section we recall the general definition for quantum affinizations and the existence of a related combinatorial fusion product introduced in \cite{her04}.

First let us define the generalized $q$-characters. For $W\in\text{Mod}(\U_q(\hat{\Hlie}))$ and an $l$-weight $(\lambda,\Psi)$ we denote:
$$W_{\lambda,\Psi}=\{v\in W_{\lambda}|\exists p\geq 0, \forall i\in I, m\geq 0, (\phi_{i,\pm m}^{\pm}-\Psi_{i,\pm m}^{\pm})^p.v=0\}.$$
As $\U_q(\hat{\Hlie})$ is commutative we have $W={\bigoplus}_{(\lambda,\Psi)\text{ $l$-weight}}W_{\lambda,\Psi}$ for all $W\in\text{Mod}(\U_q(\hat{\Hlie}))$.

\begin{defi}\label{qplu} We denote by $QP_l^+$ the set of $l$-weights $(\mu,\gamma)$ such that $\mu\in P^+-Q^+$, and there exists polynomials $Q_i(z), R_i(z)$ ($i\in I$) satisfying
$R_i(0)=Q_i(0)=1$ and in $\CC[[z]]$ (resp. in $\CC[[z^{-1}]]$): $${\sum}_{m\geq 0}\gamma_{i,\pm m}^{\pm} z^{\pm m}=
q_i^{\text{deg}(Q_i)-\text{deg}(R_i)}\frac{Q_i(zq_i^{-1})R_i(zq_i)}{Q_i(zq_i)R_i(zq_i^{-1})}.$$\end{defi}

\noindent The following result was first proved in \cite[Proposition 1]{Fre} for quantum affine algebras, and a straightforward generalization for general quantum affinizations was given in \cite[Proposition 5.4]{her04}.

\begin{prop}\label{qplplus} For $V\in\text{Mod}(\U_q(\hat{\Glie}))$ and $(\lambda,\Psi)$ an $l$-weight, $(V_{(\lambda,\Psi)}\neq \{0\})\Rightarrow
((\lambda,\Psi)\in QP_l^+)$.\end{prop}

\noindent Consider formal variables $Y_{i,a}^{\pm}$ ($i\in I, a\in\mathbb{C}^*$) and $k_{\omega}$ ($\omega\in\mathfrak{h}$). Let $A$ be the
commutative group of monomials of the form $m={\prod}_{i\in I,a\in\mathbb{C}^*}Y_{i,a}^{u_{i,a}(m)} k_{\omega(m)}$, ($u_{i,a}(m)\in\ZZ$, $\omega(m)\in \mathfrak{h}$), satisfying

1) only a finite number of $u_{i,a}(m)\in\mathbb{Z}$ are non zero,

2) for $i\in I$,
$\alpha_i(\omega(m))=r_i{\sum}_{a\in\mathbb{C}^*}u_{i,a}(m)$.

\noindent The product is given by the usual multiplication of Laurent polynomials in variables $Y_{i,a}^{\pm}$ and $k_{h+h'}=k_hk_{h'}$. This group $A$ is fixed for the rest of this paper.

\noindent For $(\mu, \gamma)\in QP_l^+$ we define $Y_{\mu,\gamma}=k_{\nu(\mu)}{\prod}_{i\in I, a\in\CC^*}Y_{i,a}^{\mu_a-\nu_a}\in A$ where
$Q_i(z)={\prod}_{a\in\CC^*}(1-za)^{\mu_a}$ and $R_i(z)={\prod}_{a\in\CC^*}(1-za)^{\nu_a}$ are given by definition \ref{qplu}.

\noindent $\Yim$ is the set of $\chi\in \ZZ^A$ such that there is a finite number of elements $\lambda_1,...,\lambda_s\in \mathfrak{h}^*$ satisfying for $m\in A$ : $$(\omega(m)\notin{\bigcup}_{j=1,\cdots,s}\nu(\mathcal{S}(\lambda_j))) \Rightarrow\text{ the coefficient of $m$ in $\chi$ is $0$}.$$ 
Note that $\Yim$ has a ring structure inherited from the group structure of $A$.

\begin{defi} The $q$-character of a module $W\in\text{Mod}(\U_q(\hat{\Hlie}))$ is defined by 
$$\chi_q(W)={\sum}_{(\lambda,\gamma)\in
QP_l^+}\text{dim}(V_{\lambda,\gamma})Y_{\lambda,\gamma}\in\Yim.$$\end{defi}

\noindent The $q$-characters were introduced in \cite{Fre} for finite dimensional representations of quantum affine algebras and the generalized $q$-characters (with the term $k_{\lambda}$) for general
quantum affinizations were constructed in \cite{her04}.

\noindent Let $\text{Rep}(\U_q(\hat{\Glie}))$ be the Grothendieck group of modules $V$ in $\text{Mod}(\U_q(\hat{\Glie}))$ which have a composition series :
\begin{equation}\label{comps} L_0=\{0\}\subset L_1\subset L_2\subset...\text{ such that }{\bigcup}_{i\geq 0}L_i=V \text{ and }\forall i\geq 0\text{ $L_{i+1}/L_i$ is simple.}\end{equation}
Note that we will prove in section \ref{decomp} that for a large class of quantum affinizations all modules in $\text{Mod}(\U_q(\hat{\Glie}))$ have such a composition series.

\noindent Let us recall some results not used for the main results of this paper, but which are related to the point of view of this paper (in
particular see section \ref{qchatens}). The following result was first proved in \cite{Fre} for quantum affine algebras and then in \cite{her04} for general quantum affinizations.

\begin{lem}$\chi_q$ gives an injective group morphism $\chi_q: \text{Rep}(\U_q(\hat{\Glie}))\rightarrow \Yim $.\end{lem}

\noindent For quantum affine algebras, the category $\text{Mod}(\U_q(\hat{\Glie}))$ has a tensor structure from the usual coproduct. It is compatible with the morphism of $q$-characters :

\begin{thm}\label{rmorph}\cite{Fre} Suppose that $\U_q(\hat{\Glie})$ is a quantum affine algebra.  Consider the ring structure on $\text{Rep}(\U_q(\hat{\Glie}))$ given by the usual
coproduct. Then $\chi_q$ is a ring morphism from $\text{Rep}(\U_q(\hat{\Glie}))$ to $\Yim$.\end{thm}

\noindent We proved in \cite{her04} that for general quantum affinizations one can observe a "combinatorial" fusion phenomena in terms of $q$-characters. More precisely we have

\begin{thm}\label{combfus}\cite{her04} For general quantum affinization $\U_q(\hat{\Glie})$, the image of $\chi_q$ is a subring of $\Yim$. The induced multiplication on $\text{Rep}(\U_q(\hat{\Glie}))$ is a fusion product (the constant structure on simple representations are positive). \end{thm}

\noindent This last result was stated in \cite{her04} when the quantized Cartan matrix $C(z)$ is invertible, but the last section of \cite{her04} implies also this result and does not use the inverse of $C(z)$. From theorem \ref{combfus}, it makes sense to define :

\begin{defi} The induced multiplication of theorem \ref{combfus} on $\text{Rep}(\U_q(\hat{\Glie}))$ is denoted by $*$ and is called combinatorial fusion product.\end{defi}

\noindent This fusion product is called combinatorial because {a priori} it makes sense only in the semi-simplified category. However in the next section a relation between the combinatorial fusion product and the tensor structure of $\text{Mod}$ is established. In particular it provides a representation theoretical interpretation of $*$.

\subsection{Generalized morphism of $q$-characters for the quantum fusion tensor category}\label{qchatens} In this section we construct a morphism of $q$-characters for the
quantum fusion tensor category $\text{Mod}$. We moreover prove that it is compatible with the monoidal structure (theorem \ref{prod}).

\noindent First let us give a $u$-deformed version of proposition \ref{qplplus}.

\begin{lem}\label{inch} Let $V\in\text{Mod}(\U_q'(\hat{\Glie}))$ and $W$ as in definition \ref{ucat}. Then $$W_{u=1}=(\CC[u^{\pm}]\otimes W)/((1-u)\CC[u^{\pm}]\otimes W)$$ 
is an object of $\text{Mod}(\U_q(\hat{\Hlie}))$ and for $(\lambda,\gamma)$ an $l$-weight, we have $(W_{u=1})_{(\lambda,\gamma)}\neq
\{0\}\Rightarrow (\lambda,\gamma)\in QP_l^+$. Moreover, for $(\lambda,\gamma)\in QP_l^+$, the dimension of $\text{dim}((W_{u=1})_{(\lambda,\gamma)})$ is
independent of the choice of $W$. In particular $\chi_q(W_{u=1})$ does not not depend of the choice of $W$.\end{lem}

\demo As $W\in\text{Mod}(\U_q(\Hlie))$ we have $W_{u=1}\in\text{Mod}(\U_q(\hat{\Hlie}))$. It follows from property v) of definition \ref{ucat} that for
$(\lambda,\gamma)$ an $l$-weight of $W_{u=1}$ there is an $(l,u)$-weight $(\lambda,\gamma(u))$ such that $(\lambda,\gamma(1))=(\lambda,\gamma)$ and $V_{(\lambda,\gamma
(u))}\neq\{0\}$. Let $v\in (V_{(\lambda,\gamma
(u))} -\{0\})$ and consider $(\tilde{V}(v))_{u=1}$. It is an integrable $\U_q(\hat{\Glie})$-module generated by $v$ (theorem \ref{aformcy}). Moreover we have $v\in
((\tilde{V}(v))_{u=1})_{(\lambda,\gamma)}$ and so it follows from proposition \ref{qplplus} that
$(\lambda,\gamma)\in QP_l^+$. For the last point we have
$$\text{dim}((W_{u=1})_{(\lambda,\gamma)})={\sum}_{\{(\lambda,\gamma(u))\text{
$(l,u)$-weight}|(\lambda,\gamma(1))=(\lambda,\gamma)\}}\text{dim}_{\CC(u)}(V_{(\lambda,\gamma (u))}).$$\qed

Because of lemma \ref{inch} it makes sense to define :

\begin{defi} We define the $q$-character of a module $V\in\text{Mod}(\U_q'(\hat{\Glie}))$ by $\chi_q(V)=\chi_q(W_{u=1})$ where $W$ is as in definition \ref{ucat} .\end{defi}

\noindent Note that $\chi_q$ can also be defined on $\text{Mod}$ by summation. We use the same notation $\chi_q$ for these morphisms of $q$-characters. We prove exactly as in lemma \ref{inch} :

\begin{lem}\label{compaform} Let $V\in\text{Mod}(\U_q'(\hat{\Glie}))$ and $\tilde{V}$ be an $\mathcal{A}$-form of $V$. Then $\chi_q(V)=\chi_q((\tilde{V})_{u=1})$.\end{lem}

\noindent The $q$-characters morphism is compatible with the tensor structure :

\begin{thm}\label{prod} We use the notations of remark \ref{toref}. For $V_1,V_2\in\text{Mod}$, we have $\chi_q(V_1\otimes_f V_2)=\chi_q(V_1)\chi_q(V_2)$.\end{thm}

\demo From definition \ref{ucat} one has subspaces $W_1\subset V_1$, $W_2\subset V_2$. We have seen in the proof of lemma \ref{proddef} that we can choose $W_1\otimes W_2$
for the subspace of $V_1\otimes_f V_2$ required in definition \ref{ucat}. So $\chi_q(V_1\otimes_f V_2)=\chi_q((W_1\otimes W_2)_{u=1})$. As a direct consequence of formulas (\ref{deltauh}), (\ref{deltauh2}), $\chi_q((W_1\otimes W_2)_{u=1})$ is equal to
$\chi_q((W_1)_{u=1})\chi_q((W_2)_{u=1})$.\qed

\noindent This result is analog to theorem \ref{rmorph} where the usual coproduct of quantum affine algebras is used. Here we use the Drinfeld coproduct, that is why this statement is almost automatic. It gives an additional clue indicating that $\otimes _f$ is an interesting tensor structure.

\noindent It follows directly from the definition of $\chi_q$ that :

\begin{lem}\label{suitee} Let $V,V'\in\text{Mod}(\U_q'(\hat{\Glie}))$ such that $V'\subset V$. Then $\chi_q(V)=\chi_q(V')+\chi_q(V/V')$.\end{lem}

\noindent Moreover :

\begin{lem}\label{projq} The image of the $q$-characters morphisms for $\text{Mod}$ and $\text{Mod}(\U_q(\hat{\Glie}))$ coincide, that is to say $\chi_q(\text{Mod})=\chi_q(\text{Mod}(\U_q(\hat{\Glie})))$. \end{lem}

\demo The inclusion $\chi_q(\text{Mod}(\U_q(\hat{\Glie})))\subset
\chi_q(\text{Mod})$ follows from the existence of the functor $i$ which allows to see $\text{Mod}(\U_q(\hat{\Glie}))$ as a subcategory of $\text{Mod}^1$ (see proposition \ref{inclu}). For the other
inclusion $\chi_q(\text{Mod})\subset\chi_q(\text{Mod}(\U_q(\hat{\Glie})))$ , we see as in the proof of lemma \ref{inch} that the
$q$-characters of cyclic modules are in
$\chi_q(\text{Mod}(\U_q(\hat{\Glie})))$ and then we can use lemma \ref{suitee}.\qed

\noindent As $\chi_q(\text{Mod}(\U_q(\hat{\Glie})))\simeq \text{Rep}(\U_q(\hat{\Glie}))$, the result of lemma \ref{projq} allows to consider a projection $p:\text{Mod}\rightarrow \text{Rep}(\U_q(\hat{\Glie}))$. It satisfies "$\chi_q=\chi_q\circ p$". Then theorem \ref{prod} implies :

\begin{cor}\label{inter} We use the notations of remark \ref{toref}. For $V_1,V_2\in\text{Mod}$, we have 
$$p(V_1\otimes_f V_2)=p(V_1)*p(V_2).$$\end{cor}

\noindent In particular the quantum fusion
tensor category $\text{Mod}$ gives a representation theoretical
interpretation of $*$.

\subsection{Further results on the category $\text{Mod}(\U_q(\hat{\Glie}))$}\label{decomp} In this section we study the existence of composition series for the modules of $\text{Mod}(\U_q(\hat{\Glie}))$. We will also study the dominant monomials (corresponding to dominant $l$-weights) occurring in $q$-characters. This study is valid for a class of Cartan matrices that we describe now : in this section (and in section \ref{finitel}) we suppose that $C$ is a generalized symmetrizable Cartan matrix satisfying the condition 
\begin{equation}\label{defc}((C_{i,j}<-1)\Rightarrow (r_i=-C_{j,i}=1))\end{equation}
where the $r_i$ define $D=\text{diag}(r_1,...,r_n)$ such that $B=DC$ is symmetric (they are fixed in section \ref{datas}). It includes quantum affine algebras and quantum toroidal algebras,
except of type $A_1^{(1)}, A_{2l}^{(2)}$ ($l\geq 0$). For $p\geq 1$,
it also includes the rank 2 matrices of the form
$\begin{pmatrix}2&-1\\-p&2\end{pmatrix}$ or
$\begin{pmatrix}2&-p\\-1&2\end{pmatrix}$ which will for example
provide completely new types of $T$-system considered in section \ref{pkr}.

\noindent First we have :

\begin{lem}\label{propc} If $C$ satisfies the property (\ref{defc}), then :

1) $\forall i,j\in I$, $r_i\geq -C_{j,i}$,

2) $(i\neq j)\Rightarrow (r_i=1\text{ or }C_{i,j}=-1\text{ or }C_{i,j}=0)$,

3) for $i\in I$ such that $r_i=1$, we have [($j\neq i$ and
$C_{i,j}\neq 0) \Rightarrow (C_{j,i}=-1)$],

4) for $i,j\in I$ such that $r_i>1$, we have [$(C_{i,j}< 0)\Rightarrow
(C_{i,j}=-1)$],

5) for $i\in I$ such that $r_i>1$ and $j\in I$ such that $C_{i,j}<0$, we have [($C_{j,i}=-1$ and $r_i=r_j$)
or ($C_{j,i}=-r_i$ and $r_j=1$)],

6) $C(z)$ is invertible.\end{lem}

\noindent Note that 2) means that $C$ is $q$-symmetrizable in the sense of \cite{her03}. In particular all properties of $q$-symmetrizable Cartan matrices established in the last section of \cite{her03} are also satisfied by $C$ (but they are not used in the present paper).

\demo 1) If $C_{j,i}\geq -1$, the result is clear because $r_i\geq 1$. For $C_{j,i}<-1$, as $B=DC$ is symmetric, we have $r_j=-C_{i,j}=1$, and so $-C_{j,i}=r_i$.

2) If $C_{i,j}\notin\{0,-1\}$ we have $C_{i,j}<-1$ and so $r_i=1$.

3) Direct consequence of 1).

4) Direct consequence of 2).

5) From 4) we have $r_i=-C_{j,i}r_j$. If $r_j>1$ we have $C_{j,i}=-1$ from 4),
and so $r_i=r_j$. If $r_j=1$, we have $C_{j,i}=-r_i$. 

6) Proved in \cite[Lemma 6.9]{her03} : we have $\text{det}(C(z))\in z^R+z^{-R}+\sum_{-R<r<R}\ZZ z^r\neq 0$ where $R=\sum_{i\in I}r_i$.\qed

Let us state the main result of this section.

\begin{thm}\label{compsp} Any module $V$ in $\text{Mod}(\U_q(\hat{\Glie}))$ has a composition series (\ref{comps}).\end{thm}

\noindent The composition series is not necessarily of finite length. Note that it implies that $\text{Rep}(\U_q(\hat{\Glie}))$ is the Grothendieck group of $\text{Mod}(\U_q(\hat{\Glie}))$, and so that $\chi_q$ encodes all information from the semi-simplified category. In this section we prove this theorem \ref{compsp}. 

\noindent Let us define for monomials the analog of the notion of dominant $l$-weight (see \cite{Fre2}).

\begin{defi} An element of $A$ is said to be dominant if the powers of the $Y_{i,a}$ are positive.\end{defi}

\noindent The following result was proved in \cite{Fre} for $\U_q(\hat{sl_2})$, in \cite{Fre2} for quantum affine algebras and in \cite{her04} in general. It provides a simplification in the study of $q$-characters because we only have to look at very particular monomials.

\begin{thm}\cite{Fre, Fre2, her04}\label{cdom} An element of $\text{Im}(\chi_q)$ is uniquely determined by the multiplicity of his dominant monomials.\end{thm}

\noindent Let $m\in A-\{1\}$. $m$ is said to be right-negative if for all $a\in\CC^*$, $j\in I$, [$(u_{j,aq^L}(m)\neq 0)\Rightarrow (u_{j,aq^L}(m)<0)$] where $L=\text{max}\{l\in\ZZ|\exists i\in I, u_{i,aq^L}(m)\neq 0\}$ (see \cite{Fre2}). A product of right-negative monomials is right-negative \cite{Fre2}.

\noindent Example : for all $i\in I, a\in\CC^*$ consider 
$$A_{i,a}^{-1}=k_{-r_i\alpha_i^{\vee}}Y_{i,aq_i^{-1}}^{-1}Y_{i,aq_i}^{-1}{\prod}_{\{(j,k)|C_{j,i}\leq -1\text{ and }k\in\{C_{j,i}+1,C_{j,i}+3,...,-C_{j,i}-1\}\}}Y_{j,aq^k}\in A.$$ 
It follows from property 1) of lemma \ref{propc} that $A_{i,a}^{-1}$ is
right-negative.

\noindent Let us define for a monomial $m$ a set $\mathcal{S}(m)$ analog to the set $\mathcal{S}(\lambda)$ (see \cite{her02} for similar definitions). 

\begin{defi} For $m\in A$, let $\mathcal{S}(m)$ be the set of monomials $m'\in A$ such
that there are $(m_0=m), m_1, ...,m_{N-1}, (m_N=m')\in A$ satisfying :

1) for all $j\in\{1,\cdots,N\}$, there is $i\in I$ such that $m_j=m_{j-1}A_{i,a_1q_i}^{-1}\cdots A_{i,a_{r_j}q_i}^{-1}$  where $r_j\geq 0$ and $a_1,\cdots,a_{r_j}\in\CC^*$,

2) for $1\leq r\leq r_j$, $u_{i,a_r}(m_{j-1})\geq |\{r'\in\{1,\cdots,r_j\}|a_{r'}=a_r\}|$ where $i,r_j$ are as in condition 1).\end{defi}

\noindent For all $m'\in \mathcal{S}(m)$, $m'm^{-1}$ is
a product of $A_{k,b}^{-1}$ (it is denoted by $m'\leq m$). Moreover $m'\in \mathcal{S}(m)$ implies $\mathcal{S}(m')\subset \mathcal{S}(m)$.

\noindent It is well-known that the weights of a $\U_q(\Glie)$-highest module of highest weight $\lambda$ are in $\mathcal{S}(\lambda)$. Let us prove the following analog refined cone property (see \cite{Fre, Fre2, Nab} for previous results).

\begin{thm}\label{cone} Let $V\in\text{Mod}(\U_q(\hat{\Glie}))$ be an $l$-highest weight module of highest monomial $m$. Then all monomials $m'$ occurring in $\chi_q(V)$ are in $\mathcal{S}(m)$ and in particular $m'\leq m$.\end{thm}

\demo A weaker statement was proved in \cite[Theorem 3.1]{her05} : for $V\in\text{Mod}(\U_q(\hat{\Glie}))$ a simple $l$-highest weight module of highest monomial $m$, all monomials $m'$ occurring in $\chi_q(V)$ satisfy $m'\leq m$. The proof was based on \cite[Lemma 3.2]{her05} and used the $\U_q(\hat{sl_2})$-Weyl modules. These modules are explicitly known and satisfy the property of theorem \ref{cone} (in the strong form). So to prove the result it suffices to rewrite word by word the proof of \cite[Theorem 3.1]{her05} where "simple-modules" is replaced by "$l$-highest weight modules" and $m_1\leq m_2$ by $m_1\in \mathcal{S}(m_2)$.\qed

\noindent In the following we denote $q^{\NN}=\{q^m|m\in\NN\}$ and $q^{\ZZ}=\{q^m|m\in\ZZ\}$.

\begin{rem}\label{qzclasse} For all $m\in A$, one can consider $(a_1,\cdots,a_r)\in(\CC^*)^r$ satisfying 

1) for $i\in I, b\in\CC^*$, $(u_{i,b}(m)\neq 0)\Rightarrow (b\in a_1q^{\NN}\cup\cdots\cup a_rq^{\NN})$,

2) $a_1,\cdots, a_r$ belong to different classes of $\CC^*/q^{\ZZ}$.\end{rem}

\begin{lem}\label{fnuh} Let $m\in A$. Let $(a_1,\cdots, a_r)$ as in remark \ref{qzclasse} and $N\in \ZZ$. Then the set 
$$\mathcal{S}_N(m)=\{m'\in \mathcal{S}(m)|\forall k\in\{1,\cdots,r\}, s\geq 0, i\in I,u_{i,a_kq^{N+s}}(m'm^{-1})=0\}$$
is finite.\end{lem}

\noindent The idea of the following
proof is that for $m'$ such that there is $N$ satisfying $m'\in\mathcal{S}_N(m)$, the minimal $N$ with this property increases when we multiply $m'$ by some monomials $A_{i,a}^{-1}$. We use a double induction on $N$ and $\alpha(m)$ defined below.

\demo First note that if $m'\in (\mathcal{S}(m)-\mathcal{S}_N(m))$, then $\mathcal{S}(m')\subset (\mathcal{S}(m)-\mathcal{S}_N(m))$. Let us prove by induction on $N$ that $\mathcal{S}_N(m)$ is finite. By definition of $\mathcal{S}(m)$, for all $N\leq 0$, $\mathcal{S}_N(m)=\{m\}$.  Consider $N > 0$. Let $M={\prod}_{i\in I, 1\leq k\leq r}Y_{i,a_k}^{u_{i,a_k}(m)}k_{u_{i,a_k}\nu(\Lambda_i)}$. First for $mM^{-1}$ we can use $(a_1q,\cdots,a_rq)$ instead of the $(a_1,\cdots,a_r)$ of remark \ref{qzclasse}. With this set of complex numbers, $\mathcal{S}_N(mM^{-1})$ becomes $\mathcal{S}_{N-1}(mM^{-1})$ which is is finite by the induction hypothesis. Let $\alpha(m)={\sum}_{\{(i,k)\in I\times \ZZ| 1\leq k\leq r, u_{i,a_k}(m)>0\}}u_{i,a_k}(m)\geq 0$. $N$ is fixed and we prove by induction on $\alpha(m)\geq 0$ that $\mathcal{S}_N(m)$ is finite. For $\alpha(m)=0$, we have $\mathcal{S}_N(m)=M\mathcal{S}_N(mM^{-1})$ which is finite by above discussion. Suppose that $\alpha(m)>0$. For $m'\in M\mathcal{S}_N(mM^{-1})$ consider 
$$\mathcal{T}_N(m')=\{m''\in (\mathcal{S}_N(m')-\{m'\})|m''(m')^{-1}\in\{A_{i,a_kq^{r_i}}^{-1}\}_{i\in I, 1\leq k\leq r}\}.$$ 
By definition $\mathcal{T}_N(m')$ is finite. For $m''\in \mathcal{T}(m')$, we have $\alpha(m'')<\alpha(m)$. So by the induction hypothesis on $\alpha$, if $m''\in\mathcal{T}_N(m')$ then $\mathcal{S}_N(m'')$ is finite. But we have $$\mathcal{S}_N(m)=(M\mathcal{S}_N(mM^{-1}))\cup {\bigcup}_{\{m''\in
\mathcal{T}(m')|m'\in M\mathcal{S}_N(mM^{-1})\}}\mathcal{S}_N(m'').$$ 
So $\mathcal{S}_N(m)$ is finite.\qed

\begin{prop}\label{finnum} Let $V\in\text{Mod}(\U_q(\hat{\Glie}))$.

1) If $V$ is an $l$-highest weight module then $\chi_q(V)$ has a finite number of dominant monomials.

2) $V$ has a finite composition series if and only if $\chi_q(V)$ has a finite number of dominant monomials.\end{prop}

\demo 1) Let $m$ be the highest weight monomial of $\chi_q(V)$, and $(a_1,\cdots,a_r)\in(\CC^*)^r$ as in remark \ref{qzclasse}. From theorem \ref{cone} all monomials of $\chi_q(V)$ are in $\mathcal{S}(m)$. Let $N\in\ZZ$ such that for all $k\in\{1,\cdots,r\}$, $s\geq 0$, $i\in I$, $u_{i,a_kq^{N+s}}(m)=0$. For $m'\in (\mathcal{S}(m)-\{m\})$, $m'm^{-1}$ is right-negative. If in addition $m'$ is dominant, for all $k\in\{1,\cdots,r\}$, $s\geq 0$, $i\in I$, $u_{i,a_kq^{N+s}}(m'm^{-1})=0$ (otherwise there would be one $b\in\CC^*,i\in I$ such that $u_{i,b}(m')=u_{i,b}(m'm^{-1}) < 0$). So a dominant monomial of $(\mathcal{S}(m)-\{m\})$ is in $\mathcal{S}_N(m)$, and we can conclude with lemma \ref{fnuh}.

2) As a simple module of the category
$\text{Mod}(\U_q(\hat{\Glie}))$ is an $l$-highest weight module
(corollary \ref{simpmod}), it follows from 1) that his
$q$-character has a finite number of dominant monomials, and we get
the same for a module with a finite composition series (lemma
\ref{suitee}). As the $q$-character of a simple module has at least one dominant monomial, the converse is clear by induction on the number
of dominant monomials with multiplicity.\qed

{\it End of the proof of theorem \ref{compsp}.} We can conclude the proof of theorem \ref{compsp} because it suffices to establish that $l$-highest weight modules have a composition series (we have proved in addition that $l$-highest weight modules have a finite composition series).\qed

\subsection{Finite length and tensor structure}\label{finitel} In this section we prove the stability of the subcategory of finite length modules for the tensor structure. As a consequence we get a tensor category in the sense of \cite{ce}.

\begin{prop}\label{fpq} Let $V_1,V_2\in\text{Mod}(\U_q(\hat{\Glie}))$ such that $\chi_q(V_1)$ and $\chi_q(V_2)$ have a finite number of dominant monomials. Then
$\chi_q(V_1)\chi_q(V_2)$ has a finite number of dominant monomials.\end{prop}

\demo It follows from the 2) of proposition \ref{finnum} that $V_1$ and $V_2$ have a finite composition series, so we can suppose that they are simple. Let
$m_1,m_2$ be the highest weight monomials respectively of $V_1$ and $V_2$. The monomials of $\chi_q(V_1)\chi_q(V_2)$ are in $\mathcal{S}(m_1)\mathcal{S}(m_2)$. Let
$(a_1,\cdots,a_r)\in(\CC^*)^r$ as in remark \ref{qzclasse} for $m_1$ {\it and} $m_2$ (one can choose $(a_1,\cdots, a_r)$ so that condition of remark \ref{qzclasse} is simultaneously satisfied for $m_1$ and $m_2$). Let $N\in\ZZ$ such that for all $k\in\{1,\cdots,r\}$, $s\geq 0$, $i\in I$, $u_{i,a_kq^{N+s}}(m_1)=u_{i,a_kq^{N+s}}(m_2)=0$. For a dominant $m'=m_1'm_2'\in (\mathcal{S}(m_1)\mathcal{S}(m_2)-\{m_1m_2\})$, $m_1'm_1^{-1}$ and
$m_2'm_2^{-1}$ are right-negative or equal to $1$, so for all $k\in\{1,\cdots,r\}$, $s\geq 0$, $i\in I$, $u_{i,a_kq^{N+s}}(m_1'm_1^{-1})=u_{i,a_kq^{N+s}}(m_2'm_2^{-1})=0$ (see the argument in the proof of 1) of proposition \ref{finnum}).
So $m_1'\in \mathcal{S}_N(m_1)$, $m_2'\in \mathcal{S}_N(m_2)$ and we can conclude with lemma \ref{fnuh}.\qed

\noindent In particular the Grothendieck ring of modules in $\text{Mod}(\U_q(\hat{\Glie}))$ with a finite composition series is stable under the fusion product $*$.

\noindent Let us prove a $u$-deformed version of proposition \ref{finnum} :

\begin{lem}\label{flen} Let $V\in\text{Mod}(\U_q'(\hat{\Glie}))$.

1) If $V$ is an $(l,u)$-highest weight module then $\chi_q(V)$ has a finite number of dominant monomials.

2) $V$ has a finite composition series if and only if $\chi_q(V)$ has
a finite number of dominant monomials.\end{lem}

\demo 1) It follows from theorem \ref{aformcy} and lemma \ref {compaform} that $\chi_q(V)$ is equal to the $q$-character of an $l$-highest weight $\U_q(\hat{\Glie})$-module and so we can use 1) of proposition \ref{finnum}.

2) We can conclude as in 1) of proposition \ref{finnum}.\qed

\noindent We can conclude with the main result of this section :

\begin{thm}\label{prectens} We use the notations of remark \ref{toref}. The subcategory of finite length modules in $\text{Mod}$ is stable under the monoidal bifunctor $\otimes_f$.\end{thm}

\demo We get the result by combining 2) of lemma \ref{flen}, theorem \ref{prod} and proposition \ref{fpq}.\qed

\section{Applications}\label{tsyst} 

In this section we develop several applications from the results of previous sections. First we prove that a simple module of the category $\text{Mod}(\U_q(\hat{\Glie}))$ is a quotient of fusion of fundamental representations (proposition \ref{mainapp}). This is analog to a result of Chari-Pressley for the usual tensor product of quantum affine algebras (see \cite[Corollary 12.2.8]{Cha2}). We establish a cyclicity property ($l,u$-highest weight property of fusion of $l$-highest modules) which allows to control the "size" of this fusion module (theorem \ref{lhw} and corollary \ref{conts}) : this property is completely different compared to the usual tensor product of quantum affine algebras (see remark \ref{usucop}), and the corresponding proof is independent of it. Then we establish exact sequences involving fusion of Kirillov-Reshetikhin modules
(theorem \ref{exact}). It is related to new generalized $T$-systems (theorem \ref{gentsyst}). For quantum affine algebras, such $T$-systems and related exact sequences involving usual tensor product of Kirillov-Reshetikhin modules were first established in \cite{Nad} for simply-laced cases and then in \cite{her06} for the general case with a different proof. The proof of the new $T$-system is a slight modification of the proof of \cite{her06}, but the proof of the corresponding exact sequence uses the fusion tensor product of the present paper. In subsection \ref{furth} we address complements and question that we hope to discuss in another paper.

In this entire section we use the notations of remark \ref{toref}.

\subsection{Application to $l$-highest weight modules} In this section we first prove a general result analog to \cite[Corollary 12.2.8]{Cha2}, and then we prove the following new cyclicity property which has no analog for the usual tensor product : for $V$, $V'$ $l$-highest weight modules, $i(V)\otimes_f i(V')$ is $(l,u)$-highest weight (theorem \ref{lhw}).

\subsubsection{Statement of the main results} First we have :

\begin{prop}\label{mainapp} A simple module $V$ of the category $\text{Mod}(\U_q(\hat{\Glie}))$ is a quotient of a fusion of fundamental representations $L_{\mu}*_f L_{i_1,a_1}*_f\cdots *_f L_{i_M,a_m}$. Moreover we have 
$$\chi_q(L_{\mu}){\prod}_{m=1,\cdots,M}\chi_q(L_{i_1,a_1})=\chi_q(L(\lambda,\Psi))+E$$
where $E\in\Yim$ has positive coefficients.\end{prop}

\demo From corollary \ref{simpmod}, $V$ is of the form $L(\lambda,\Psi)$ where $(\lambda,\Psi)$ is a dominant $l$-weight. Let $(P_i(u))_{i\in I}$ be the corresponding set of Drinfeld polynomials. Consider $((i_1,a_1),\cdots,(i_M,a_M))\in (I\times\CC^*)^M$ such that for each $i\in I$, $P_i(u)={\prod}_{\{m|i_m=i\}}(1-ua_m)$. From corollary \ref{infosf}, the fusion module
$$L_{(\lambda-{\sum}_{i\in I}\lambda(\alpha_i^{\vee})\Lambda_i)}*_f L_{i_1,a_1}*_f L_{i_2,a_2}*_f \cdots *_f L_{i_M,a_M}$$
is an $l$-highest weight module of $l$-highest weight $(\lambda,\Psi)$ and so admits $L(\lambda,\Psi)$ as a quotient.

\noindent From the first part and by lemma \ref{suitee}, we have $\chi_q(L_{\mu}*_f L_{i_1,a_1}*_f\cdots *_f L_{i_M,a_m})=E'+\chi_q(L(\lambda,\Psi))$ where $E'$ has positive coefficients. But $L_{\mu}*_f L_{i_1,a_1}*_f\cdots *_f L_{i_M,a_m}$ comes from an $\mathcal{A}$-form of an $\U_q'(\hat{\Glie})$-submodule of
$L_{\mu}\otimes_f L_{i_1,a_1}\otimes_f L_{i_2,a_2}\otimes_f \cdots\otimes_f L_{i_M,a_M}$. So we can conclude with theorem \ref{prod}. \qed

\noindent The first part of proposition \ref{mainapp} is analog to a theorem of Chari-Pressley
related to the usual tensor product for quantum affine algebras (see
for example \cite[cor. 12.2.8]{Cha2}).

\noindent We can control the "size" of the fusion of $l$-highest weight modules. Indeed let us state the main result of this section which is a cyclicity property.

\begin{thm}\label{lhw} Let $V,V'\in\text{Mod}(\U_q(\hat{\Glie}))$ be two $l$-highest weight modules. Then $i(V)\otimes_f i(V')$ is an $(l,u)$-highest weight
module.\end{thm}

\noindent This theorem is proved in section \ref{lhwm}. This
result is
completely different from the properties of the usual tensor product of
quantum affine algebras (see remark \ref{usucop}, \cite[Theorem 4]{c} and \cite[Theorem 9.1]{kas}). Here all fusion
products are $l$-highest weight. Note that evidences of this result were observed in the examples of section \ref{exkr}.

\noindent We have the following direct consequence of Theorem \ref{lhw}.

\begin{cor}\label{conts} Let $V_1,V_2,\cdots,V_r\in\text{Mod}(\U_q(\hat{\Glie}))$ be $l$-highest weight modules. Then the $l$-highest weight module $V=V_1*_f V_2*_f\cdots*_f V_r$ satisfies $\chi_q(V)=\chi_q(V_1)\chi_q(V_2)\cdots\chi_q(V_r)$.\end{cor}

\noindent In general (for example for quantum toroidal algebras)
the existence of such "big" $l$-highest weight modules was not known so far. 

\noindent For quantum affine
algebras, Weyl modules are certain maximal finite dimensional $l$-highest weight 
modules for fixed Drinfeld polynomial. Corollary \ref{conts} implies that the dimension of a Weyl module is
larger that the product of the dimension of the corresponding fundamental representations (this result is already known, consequence of \cite[Theorem 4]{c} and consequence of \cite[Theorem 9.1]{kas}). By analogy, for general quantum
affinizations, a maximal integrable $l$-highest weight module would
have such a fusion of fundamental representations as a quotient. In a paper in preparation
these fusion modules will be used to study the block decomposition of
$\text{Mod}(\U_q(\hat{\Glie}))$ (see section \ref{odev}). They
are also used in the second application (in section \ref{applikr}).

\subsubsection{Proof of theorem \ref{lhw}}\label{lhwm} We first prove preliminary lemmas. Let us recall the following well-known result.

\begin{lem}\label{linind} Let $M\geq 0$. The sequences $((\lambda^m m^k)_{m\geq M})_{\lambda\in\CC^*, k\geq 0}$ are $\CC$-linearly independent. \end{lem}

For convenience of the reader, we give a proof.

\demo For $k\geq 0$ consider $G_k(z)={\sum}_{m\geq 0}z^mm^k$. We
have $G_0(z)=1/(1-z)$ and $G_{k+1}(z)=zG_k'(z)$. So by induction on $k$, we have $G_k(z)=P_k(z)/(1-z)^{k+1}$ where $P_k(z)\in\CC[z]$ and $P_k(1)\neq
0$. So $G_k(z)\in 1/(1-z)^{k+1}+\sum_{k'<k}\CC.G_{k'}(z)$ and $(G_k(\lambda z))_{\lambda\in\CC^*, k\geq 0}$ is a
$\CC$-basis of $\CC(z)/\CC[z^{\pm}]$. $M\geq 0$ is fixed and consider $G_k^{(M)}={\sum}_{m\geq M}z^mm^k$. We have
$G_k^{(M)}(z\lambda)\in G_k(z\lambda)+\CC[z^{\pm}]$ and so the
$G_k^{(M)}(z\lambda)$ are linearly independent in
$\CC(z)/\CC[z^{\pm}]$.\qed

\noindent As a consequence, we have :

\begin{lem}\label{annclas} Let $V$ be a $\CC(u)$-vector space generated by vectors $(f_{\lambda,k,k'})_{\lambda\in\CC^*,k\geq 0, k'\geq 0}$.
We suppose that only a finite number of these vectors are not equal
to $0$. For $m\geq 0$ we consider
$\mu_m={\sum}_{\lambda\in\CC^*,k\geq 0,k'\geq 0}\lambda^m
m^k u^{-mk'}f_{\lambda,k,k'}\in V$. Then $V$ is generated by the
$(\mu_m)_{m\geq 0}$.\end{lem}

\demo As $V$ is a finite dimensional vector space we can prove the
result by induction on $\text{dim}(V)$. For $V=\{0\}$ it is clear.
In general let $V'\subset V$ be the subspace generated by the $\mu_m$, $m\geq 0$.
Suppose that there is $f_{\lambda_0,k_0,k_0'}\notin V'$. By 
the induction hypothesis the image of $V'$ in $V/(\CC(u)f_{\lambda_0,k_0,k_0'})$ generates the
space $V/(\CC(u)f_{\lambda_0,k_0,k_0'})$. So $\text{dim}(V') = \text{dim}(V)-1$. Consider
$\phi:V\rightarrow \CC(u)$ a non zero linear map such that $\text{Ker}(\phi)=V'$. We have
for all $m\geq 0$, ${\sum}_{\lambda\in\CC^*,k\geq 0,k'\geq
0}\lambda^m m^k u^{-mk'}\phi(f_{\lambda,k,k'})=0$.
As there is a finite number of non zero $f_{\lambda,k,k'}$, we can consider $P(u)\in(\CC(u)-\{0\})$, such that $\phi(f_{\lambda,k,k'})P(u)=P_{\lambda,k,k'}(u^{-1})$ where
$P_{\lambda,k,k'}\in\CC[u]$. As $\phi\neq 0$, there is at least one polynomial
$P_{\lambda,k,k'}\neq 0$. Consider a maximal $k'=K'$ such that there is $(\lambda,k)$ satisfying $P_{\lambda,k,K'}\neq 0$. Let $\Lambda$ be the set of $(\lambda,k,K')$ such that the polynomials
$P_{\lambda,k,K'}$ have a maximal degree. Let $\alpha_{\lambda,k,K'}$ be the leading coefficient of
$P_{\lambda,k,K'}$. There is $M\geq 0$ such
that for $m\geq M$, $0={\sum}_{(\lambda ,k,
K')\in\Lambda}\lambda^m m^k\alpha_{\lambda,k,K'}$. Contradiction because from lemma \ref{linind}, the sequences
$(\lambda^m m^k)_{m\geq M}$ are $\CC$-linearly independent.\qed

\noindent In the following we denote $\U_q'(\hat{\Hlie})=\U_q(\hat{\Hlie})\otimes \CC(u)\subset \U_q'(\hat{\Glie})$.

\begin{lem}\label{aidecyc} Let $W$ be an object of $\text{Mod}(\U_q(\hat{\Hlie}))$ and $W' \subset (W\otimes \CC(u))$ a $\CC(u)$-subspace of $W\otimes \CC(u)$. We consider $\mathcal{R}\in \U_q'(\hat{\Hlie})\otimes\CC(u)$ such that $\mathcal{R}(W\otimes \CC(u))\subset W'$ and which can be expanded in the form $\mathcal{R}=\sum_{r\geq 0} u^r g_r$ where $g_0\in \CC^*$ and for $r\geq 1$, $g_r\in\U_q(\hat{\Hlie})$. Then $W'=W\otimes\CC(u)$.\end{lem}

\demo It suffices to prove that for $w$ an $l$-weight vector of $W$, we have $w\in W'$. The restriction of $\mathcal{R}$ on $\U_q'(\hat{\Hlie}).w$ is of the form $(\gamma \text{Id}+ n)$ where $\gamma\in\CC(u)$ and $n$ is nilpotent. Consider $\gamma_r\in\CC$ the
eigenvalue of $g_r$ on $\U_q'(\hat{\Hlie}).w$. As the $g_r$ commute, we have
$\gamma=g_0+{\sum}_{r>0}u^r\gamma_r$. As $g_0\neq 0$, $\gamma$ is invertible in $\CC(u)$. So $w'=\gamma^{-1}{\sum_{v\geq 0}}(n/\gamma)^v.w$ makes sense
in $\U_q'(\hat{\Hlie}).w$. We can conclude as $\mathcal{R}.w'=w$ and by hypothesis $\mathcal{R}.w'\in W'$.\qed

{\it End of the proof of theorem \ref{lhw} :}

Consider $v\in V_{\lambda}$ (resp. $v'\in V_{\lambda'}$) an
$l$-highest weight vector of $V$ (resp. of $V'$). Let $W=\U_q'(\hat{\Glie}).(v\otimes v')$. We prove that $W=i(V)\otimes_f i(V')$. For $\mu$ a weight of
$i(V)\otimes_f i(V')$, we have $\mu\in\lambda+\lambda'-Q^+$. The map $h$ is defined in section \ref{datas}. Let us
prove by induction on $h(\lambda+\lambda'-\mu)\geq 0$ that
$(i(V)\otimes_f i(V'))_{\mu}\subset W$. For $h(\lambda+\lambda'-\mu)=0$ it is clear. In general let $\mu_1,\mu_2$ such that $\mu=\mu_1+\mu_2$ and $(i(V)_{\mu_1}\otimes_f i(V')_{\mu_2})\subset W$. It
suffices to prove that ${\sum}_{i\in I,m\in\ZZ}(x_{i,m}^-(V_{\mu_1})\otimes V_{\mu_2}')\subset W$ and ${\sum}_{{i\in I, m\in\ZZ}}(V_{\mu_1}\otimes
x_{i,m}^-(V_{\mu_2}'))\subset W$.

\noindent Consider $w\in V_{\mu_1}$, $w'\in
V_{\mu_2}'$. Let $p\in\ZZ$. Formulas (\ref{deltaup}) and
(\ref{defit}) give for $m\in\ZZ$ : $x_{i,p-m}^-.(w\otimes w')=A+B+C$
where $A=u^{p-m}w \otimes (x_{i,p-m}^-.w')$ , $B=(x_{i,p-m}^-.w)\otimes
(k_i.w')$ and $C={\sum}_{r>0}u^r((x_{i,p-m-r}^-.w)\otimes
(\phi_{i,r}^+.w'))$.

\noindent From proposition \ref{inclu}, the properties of definition \ref{ucat} are satisfied by $i(V)$ and $i(V')$. So let us consider the decompositions of property iv) of definition \ref{ucat} : we have for all $m\geq 0$, all $w\in V_{\mu_1}$, all $w'\in V_{\mu_2}'$
$$x_{i,p-m}^-(w)={\sum}_{\lambda\in\CC^*, k\geq 0}\lambda^m
m^k f_{\lambda,k}(w),$$
$$x_{i,p-m}^-(w')={\sum}_{\lambda\in\CC^*, k\geq
0}\lambda^m m^k f_{\lambda,k}'(w').$$ 
Moreover we can consider the decomposition of lemma \ref{phim} : for $r\geq 1$,
$\phi_{i,r}^+(w')={\sum}_{\lambda\in\CC^*, k\geq 1}\lambda^r
r^k g_{\lambda,k}(w')$. Only a finite number of $f_{\lambda,k},
f_{\lambda,k}', g_{\lambda,k}$ are not equal to zero. 

\noindent We have :
$$A={\sum}_{\lambda\in\CC^*, k\geq
0}u^{-m}\lambda^m m^k (u^p w\otimes f_{\lambda,k}'(w')),$$
$$B={\sum}_{\lambda\in\CC^*, k\geq 0}\lambda^m m^k
(f_{\lambda,k}(w)\otimes q_i^{\mu_2(\alpha_i^{\vee})}w'),$$
$$C={\sum}_{\lambda\in\CC^*, k\geq
0,\lambda'\in\CC^*, k'\geq 0, r>0}u^r\lambda^{m+r}(m+r)^k
(\lambda')^r r^{k'} f_{\lambda,k}(w)\otimes g_{\lambda',k'}(w')$$
$$={\sum}_{\lambda\in\CC^*, k\geq
0, s=0,\cdots,k} \lambda^m m^{k-s} f_{\lambda,k}(w)\otimes
\mathcal{R}_{s,\lambda}(w'),$$ 
where 
$$\mathcal{R}_{s,\lambda}=\begin{bmatrix}k\\s\end{bmatrix}{\sum}_{\lambda'\in\CC^*, k'\geq 0,
r>0}u^r\lambda^rr^s (\lambda')^r r^{k'} g_{\lambda',k'}.$$
Note that $\mathcal{R}_{s,\lambda}$ is independent of $m$. As $\phi_{i,r}^+$ equals ${\sum}_{\lambda\in\CC^*, k\geq 1}\lambda^rr^k g_{\lambda,k}$ in $\text{End}(V_{\mu_2}')$, it follows from lemma \ref{annclas} that $g_{\lambda,k}$ viewed in $\text{End}(V_{\mu_2}')$ is an operator of $\U_q(\hat{\Hlie})$. And so $\mathcal{R}_{s,\lambda}$ is also an operator of $\U_q(\hat{\Hlie})$ viewed in $\text{End}(V_{\mu_2}')$.

\noindent As $A+B+C\in W$, it follows from lemma \ref{annclas} that for all
$\lambda\in\CC^*, K\geq 0$ :
\begin{equation}\label{premw}V_{\mu_1}\otimes f_{\lambda,K}'(V_{\mu_2}')\subset W,\end{equation}
\begin{equation}\label{deuxw}(f_{\lambda,K}\otimes
q_i^{\mu_2(\alpha_i^{\vee})}\text{Id}+{\sum}_{\{(k,s)|0\leq s, k-s=K\}}f_{\lambda,k}\otimes
\mathcal{R}_{s,\lambda})(V_{\mu_1}\otimes V_{\mu_2}')\subset W.\end{equation}

\noindent From inclusion (\ref{premw}) we get ${\sum}_{i\in I,
m\in\ZZ}(V_{\mu_1}\otimes x_{i,m}^-(V_{\mu_2}'))\subset W$.

\noindent Let us prove that ${\sum}_{i\in I,
m\in\ZZ}(x_{i,m}^-(V_{\mu_1})\otimes V_{\mu_2}')\subset W$. We fix $\lambda\in\CC^*$. Let $K_0$ be the maximal $k$ such that $f_{\lambda,k}\neq 0$. We prove by induction on $K_0-K\geq 0$ that
$(f_{\lambda,K}(V_{\mu_1})\otimes V_{\mu_2}')\subset W$. Let $w'$ be an $l$-weight vector in $V_{\mu_2}'$. For
$K=K_0$, inclusion (\ref{deuxw}) gives $f_{\lambda,K_0}(V_{\mu_1})\otimes
(q_i^{\mu_2(\alpha_i^{\vee})}w'+\mathcal{R}_{0,\lambda}(w'))\subset W$. From the defining formula of $\mathcal{R}_{0,\lambda}$, we can use lemma \ref{aidecyc} and we get $(f_{\lambda,K_0}(V_{\mu_1}')\otimes w')\subset W$. So $(f_{\lambda,K_0}(V_{\mu_1}')\otimes V_{\mu_2}')\subset W$. 
For $K\leq K_0$, $w\in V_{\mu_1}, w'\in V_{\mu_2}'$, the vector $f_{\lambda,K}(w)\otimes
((q_i^{\mu_2(\alpha_i^{\vee})}+\mathcal{R}_{0,\lambda}).w')$ is equal
to 
$$(f_{\lambda,K}(w)\otimes
q_i^{\mu_2(\alpha_i^{\vee})}w'+{\sum}_{\{(k,s)|k-s=K\}}f_{\lambda,k}(w)\otimes
\mathcal{R}_{s,\lambda}(w'))-{\sum}_{\{(k,s)|k-s=K,k>K\}}f_{\lambda,k}(w)\otimes
\mathcal{R}_{s,\lambda}(w').$$ 
By inclusion (\ref{deuxw}) the first term in in $W$. The induction hypothesis gives that the second term is in $W$. So $f_{\lambda,K}(w)\otimes ((q_i^{\mu_2(\alpha_i^{\vee})}+\mathcal{R}_{0,\lambda}).w')\in W$. If moreover $w'$ is an $l$-weight vector, from the defining formula of $\mathcal{R}_{0,\lambda}$ we can use lemma \ref{aidecyc} and we get $(f_{\lambda,K}(w)\otimes w')\in W$. So $(f_{\lambda,K}(V_{\mu_1})\otimes V_{\mu_2}')\subset W$ and ${\sum}_{i\in I, m\in\ZZ}(x_{i,m}^-(V_{\mu_1})\otimes V_{\mu_2}')\subset W$.\qed

%But $q_i^{\mu_2(\alpha_i^{\vee})}w'+\mathcal{R}_{0,\lambda}(w')=\phi.w'$
%where $\phi\in\U_q'(\hat{\Hlie})$ (as we observed that $\mathcal{R}_{0,\lambda}$ is an operator of $\U_q(\hat{\Hlie})$ viewed in $\text{End}(V_{\mu_2}')$). As $w'$ is an $l$-weight vector,
%the action of $\phi$ restricted to $\U_q'(\hat{\Hlie}).w'$ is of the form
%$(\gamma \text{Id}+ n)$ where $\gamma\in\CC(u)$ and $n$ is nilpotent. Consider $\gamma_{\lambda',k'}\in\CC$ the
%eigenvalue of $g_{\lambda',k'}$. As the $g_{\lambda',k'}$ commute as operators on $V_{\mu_2'}$, we have
%$\gamma=q_i^{\mu_2(\lambda_i^{\vee})}+{\sum}_{\lambda'\in\CC^*,k'\geq
%0}\gamma_{\lambda',k'}{\sum}_{r>0}(u\lambda\lambda')^r
%r^{k'}$ which is invertible in $\CC(u)$. So for $H\in\U_q'(\hat{\Hlie})$,
%$w''=\gamma^{-1}{\sum_{v\geq 0}}(n/\gamma)^v.Hw'$ makes sense
%in $\U_q'(\hat{\Hlie}).w'$. Relation (\ref{deuxw}) gives
%$$f_{\lambda,K_0}(w)\otimes (q_i^{\mu_2(\alpha_i^{\vee})}w''+\mathcal{R}_{0,\lambda}(w''))\in W.$$
%But this term is equal to $f_{\lambda,K_0}(w)\otimes \phi.w'' = f_{\lambda,K_0}(w)\otimes H.w'$. As this property holds for all $H\in\U_q'(\hat{\Hlie})$, 

\subsection{Application to exact sequences involving generalized Kirillov-Reshetikhin modules}\label{applikr} Exact sequences involving usual tensor products of Kirillov-Reshetikhin modules and related to $T$-systems were first established in \cite{Nab, Nad} for simply-laced quantum affine algebras and in \cite{her06} with a different proof for general quantum affine algebras. These relations are of particular importance because they imply the Kirillov-Reshetikhin conjecture which predicts explicit formulas for the characters of the Kirillov-Reshetikhin modules (see \cite{Nad, her06}). In this section we establish analog exact sequences for a large class of general quantum affinizations : we replace the usual tensor product of quantum affine algebras by the fusion tensor product $*_f$ introduced in this paper for quantum affinizations and we use the cyclicity property established in theorem \ref{lhw}.

\subsubsection{Statement of the main result}  We suppose that $C$ is a generalized symmetrizable Cartan matrix satisfying the condition (\ref{defc}) (it includes quantum affine
algebras and quantum toroidal algebras, except for $A_1^{(1)}, A_{2l}^{(2)}$, $l\geq 0$).

\noindent For $r\geq 0, a\in\CC^*, i\in I$ consider the following fusion module  \begin{equation}\label{skai}S_{r,a}^{(i)}=L_{\nu_{r,a}^{(i)}}*_f({*_f}_{\{(j,k)|C_{j,i}<0\text{ , }1\leq k\leq -C_{i,j}\}}W_{-C_{j,i}+E(r_i(r-k)/r_j),aq_j^{-(2k-1)/C_{i,j}}}^{(j)})\end{equation}
where
$$\nu_{r,a}^{(i)}=r(\Lambda_i-\alpha_i)-{\sum}_{\{(j,k)|C_{j,i}<0\text{ , }1\leq k\leq -C_{i,j}\}}(-C_{j,i}+E(r_i(r-k)/r_j))\Lambda_j.$$
The motivation for these formulas is that lemma \ref{motivform} studied below is satisfied. A priori, the $\U_q(\hat{\Glie})$-module $S_{r,a}^{(i)}$ is not
well-defined because it could depend of the order of the terms, but :

\begin{thm}\label{exact} Let $i\in I, a\in\CC^*, r\geq 1$. Then

(1) the module $S_{r,a}^{(i)}$ is well-defined and simple,

(2) the module $W_{r+1,a}^{(i)}*_f W_{r-1,aq_i^2}^{(i)}$ is simple,

(3) there exists an exact sequence 
$$0\rightarrow S_{r,a}^{(i)}\rightarrow W_{r,a}^{(i)}*_f W_{r,aq_i^2}^{(i)}\rightarrow W_{r+1,a}^{(i)}*_f
W_{r-1,aq_i^2}^{(i)}\rightarrow 0.$$\end{thm}

\noindent These results appeared for quantum affine algebras, with $\otimes$ instead of $*_f$, first in \cite{Nad} (simply-laced cases) and then in \cite{her06} (general cases).

\noindent As $W_{k,a}^{(i)}*_f W_{k,aq_i^2}^{(i)}$ is $l$-highest weight (corollary \ref{infosf}), theorem \ref{exact} implies that it is not semi-simple.

\noindent Theorem \ref{exact} is proved in this section \ref{applikr}. The
proof of the corresponding result for $q$-characters is
a slight modification of the proof of \cite{her06}, but the proof of the exact sequence uses the new developments
of the present paper.

\subsubsection{Preliminary results}

\begin{defi} A $\U_q(\hat{\Glie})$-module of $\text{Mod}(\U_q(\hat{\Glie}))$ is said to be special if his $q$-character has a unique dominant monomial.\end{defi}

\begin{cor}\label{spefus} Let $\{(\lambda_1,\Psi_1),\cdots,(\lambda_M,\Psi_M)\}$ be a set of dominant $l$-weights and $(\lambda,\Psi)$ such that $Y_{\lambda,\Psi}={\prod}_{m=1,\cdots,M}Y_{(\lambda_m,\Psi_m)}$. Suppose that only
the unique dominant monomial $Y_{\lambda,\Psi}$ appears in ${\prod}_{m=1,\cdots,M}\chi_q(L(\lambda_m,\Psi_m))$. Then we have
$$L(\lambda,\Psi)\simeq L(\lambda_{\sigma(1)},\Psi_{\sigma(1)})*_f L(\lambda_{\sigma(2)},\Psi_{\sigma(2)})*_f
\cdots*_fL(\lambda_{\sigma(M)},\Psi_{\sigma(M)}).$$
In particular this fusion module is special and independent of the
order of the terms.
\end{cor}

\demo In the second statement of proposition \ref{mainapp} we have
necessarily $E=0$ because of theorem \ref{cdom}.\qed

\subsubsection{New generalized $T$-systems}\label{pkr} The proof of theorem
\ref{exact} is based on the fact that the $q$-characters of
Kirillov-Reshetikhin modules solve the following generalized
$T$-system :

\begin{thm}\label{gentsyst} Let $a\in\CC^*, r\geq 1, i\in I$. We have 
$$\chi_q(W_{r,a}^{(i)})\chi_q(W_{r,aq_i^2}^{(i)})=\chi_q(W_{r+1,a}^{(i)})\chi_q(W_{r-1,aq_i^2}^{(i)})+\chi_q(S_{r,a}^{(i)}),$$
$$\text{ where }\chi_q(S_{r,a}^{(i)})=\chi_q(L_{\nu_{r,a}^{(i)}}){\prod}_{\{(j,k)|C_{j,i}<0\text{ , }1\leq k\leq -C_{i,j}\}}\chi_q(W_{-C_{j,i}
+E(r_i(r-k)/r_j),aq_j^{-(2k-1)/C_{i,j}}}^{(j)}).$$
\end{thm}

\noindent (Note that this result is also satisfied for the algebra $\U_q(\hat{\Glie})$ with affine quantum Serre relations (\ref{equadeuxc}) as it only involves only $q$-characters and so makes sense in the Grothendieck group. In particular for this result the restriction of remark \ref{toref} is not necessary.)

\noindent For quantum affine algebras $T$-systems were first established in \cite{Nad} (simply-laced cases) and for all types in \cite{her06} with a different proof. They can be considered as "induction" systems of relations for the $q$-characters of Kirillov-Reshetikhin modules. In this paper we follow the proof of \cite{her06} (whose parts of the plan first appeared in \cite{Nad}) to get the
generalized result of theorem \ref{gentsyst}. Note that in the
particular case of simply-laced quantum toroidal algebras, geometric constructions
are discussed in \cite[Section 6]{Nac}; so in analogy to the case of simply-laced quantum affine algebras, there could be an alternative geometric proof for the $T$-systems of simply-laced quantum toroidal algebras.

\noindent This definition of a $T$-system for this larger class of Cartan matrix is new. It coincides with the definition for quantum affine algebras \cite{kns}. Indeed we have :

\begin{lem}\label{coint} For $i\in I$ such that $r_i>1$, we have 
$$\chi_q(S_{r,a}^{(i)})=\chi_q(L_{\nu_{r,a}^{(i)}})({\prod}_{\{j\in I|C_{j,i}=-1\}}\chi_q(W_{r,aq_i}^{(j)}))({\prod}_{\{j\in I|C_{j,i}=-r_i\}}\chi_q(W_{r_ir,aq}^{(j)})).$$
For $i\in I$ such that $r_i=1$, we have 
$$\chi_q(S_{r,a}^{(i)})=\chi_q(L_{\nu_{r,a}^{(i)}})({\prod}_{\{(j,k)|j\in I,C_{j,i}=-1\text{ and }k\in\{1,\cdots,r_j\}\}}\chi_q(W_{1+E((r-k)/r_j),aq^{2k-1}}^{(j)})).$$
\end{lem}

\demo Let $i \in I,
a\in\CC^*, r\geq 0$ fixed. For $j$ such that $C_{j,i}<0$ denote 
$$X_j={\prod}_{1\leq k\leq -C_{i,j}}\chi_q(W_{-C_{j,i}
+E(r_i(r-k)/r_j),aq_j^{-(2k-1)/C_{i,j}}}^{(j)}).$$
We have $\chi_q(S_{r,a}^{(i)})=\prod_{\{j\in I|C_{j,i}<0\}}X_j$. 

\noindent First suppose that $r_i>1$. Property 4) of lemma \ref{propc} gives $C_{i,j}=-1$ and so $X_j=\chi_q(W_{1+E(r_i(r-1)/r_j),aq^{r_j}}^{(j)})$. Property 5) of lemma \ref{propc} gives that ($C_{j,i}=-1$ or $C_{j,i}=-r_i$). 

If $C_{j,i}=-1$, we have $r_i=r_j$ and 
$$X_j=\chi_q(W_{1+E(r_i(r-1)/r_j),aq_j^{-1/C_{i,j}}}^{(j)})=\chi_q(W_{1+E(r-1),aq^{r_j}}^{(j)})=\chi_q(W_{r,aq_i}^{(j)}).$$

If $C_{j,i}=-r_i$, we have $r_j=1$ and 
$$X_j=\chi_q(W_{-C_{j,i}
+E(r_i(r-1)/r_j),aq_j^{-1/C_{i,j}}}^{(j)})=\chi_q(W_{r_i+E(r_i(r-1)),aq}^{(j)})=\chi_q(W_{r_ir,aq}^{(j)}).$$

\noindent Secondly suppose that $r_i=1$. From property 3) of lemma \ref{propc}, we have $C_{j,i}=-1$, $r_j=-C_{i,j}$ and 
$$X_j={\prod}_{k=1,\cdots,r_j}\chi_q(W_{1+E((r-k)/r_j),aq^{-r_j(2k-1)/C_{i,j}}}^{(j)})={\prod}_{k=1,\cdots,r_j}\chi_q(W_{1+E((r-k)/r_j),aq^{2k-1}}^{(j)}).$$
\qed

\subsubsection{End of the proof of theorem \ref{gentsyst}}\label{spfun}

We use a slight modification of the proof of \cite[Theorem 3.4 (1)]{her06}. A part of the plan of this proof (in particular theorem \ref{krspe}) first appeared in \cite{Nad} for simply-laced quantum affine algebras with different geometric arguments inside the proof which can be used only for simply-laced cases so far. In this section we recall the main steps of this proof. Moreover we give the list of results of \cite{her06} whose proof has to be slightly modified for the general situation considered in theorem \ref{gentsyst}.

Let us denote by $m_{k,a}^{(i)}$ the highest weight monomial of $\chi_q(W_{k,a}^{(i)})$. We have the following three steps :

1) The first step is an analog of \cite[Lemma 4.4]{her06}.

\begin{lem}\label{spec} For $m$ a monomial occurring in $\chi_q(W_{k,a}^{(i)})$, $((m\neq m_{k,a}^{(i)})\Rightarrow (m\leq mA_{i,aq_i^{2k-1}}^{-1}))$.\end{lem}

\noindent In the proof of this result, first an analog of the technical result \cite[Lemma 4.3]{her06} is needed. But it follows from lemma \ref{propc} that the quantized Cartan matrix is invertible, and so we can use \cite[Lemma 5.10]{her04} in this situation. To prove lemma \ref{spec}, we also need

\begin{lem} Let $\lambda$ be a weight and $j\in I$ such that $(W_{k,a}^{(i)})_{\lambda}\neq \{0\}$ and $\lambda=k\Lambda_i-\alpha_j$. Then $i=j$.\end{lem}

\demo It follows from \cite[Corollary 3.10]{her05} that the result is true for $k=1$. So the result follows from proposition \ref{mainapp}.\qed

\noindent Then we can conclude word by word the proof of lemma \ref{spec} as in \cite[Lemma 4.4]{her06}.

2) As a direct consequence of lemma \ref{spec}, we have :

\begin{thm}\label{krspe} The Kirillov-Reshetikhin modules are special.\end{thm}

\noindent This result is useful for our study because from theorem \ref{cdom}, we get a characterization of the $q$-character of Kirillov-Reshetikhin modules. 

3) The last step is to prove the $T$-system by using theorem \ref{cdom}. So we need to determine the dominant monomials of the $q$-characters involved in the $T$-system. First we have an analog of the first part of \cite[Proposition 5.3]{her06}.

\begin{lem}\label{sspe} $S_{r,a}^{(i)}$ is special.\end{lem}

\demo To use the proof of \cite{her06}, we need an analog of \cite[Lemma 5.2]{her06} which is a technical combinatorial lemma on monomials of $\chi_q(S_{r,a}^{(i)})$. We use the notations of \cite{her06}, except that we denote by $r$ the $k$ of \cite{her06}. The inequality $\beta=\mu_a(A_{j_l,aq_{j_l}^{2k_l-1}}^{-1})>\mu_a(\alpha)$ has to be studied. First we have $\mu_a(\alpha)\leq 2r_ir-1$. Then by using the description of lemma \ref{coint}, we get

if $r_i>1$ and $r_j=r_i$, then $\beta=r_i+2r_jr=r_i+2r_ir$,

if $r_i>1$ and $r_j=1$, then $\beta=1+2r_ir_jr=1+2r_ir$,

if $r_i=1$, then $\beta=2k-1+2r_j(1+E((r-k)/r_j))> 2k-1+2r_j((r-k)/r_j)=2rr_i-1$.\qed

\noindent One can check exactly as in \cite{her06} that $\chi_q(W_{r,a}^{(i)})\chi_q(W_{r,aq_i^2}^{(i)})-\chi_q(W_{r+1,a}^{(i)})\chi_q(W_{r-1,aq_i^2}^{(i)})$ has a unique dominant monomial.

\noindent To conclude the proof of the $T$-system, we need an analog of the end of the statement of \cite[Proposition 5.3]{her06}, that is to say :

\begin{lem}\label{motivform} The highest weight monomial of $\chi_q(S_{r,a}^{(i)})$ is the unique dominant monomial of $\chi_q(W_{r,a}^{(i)})\chi_q(W_{r,aq_i^2}^{(i)})-\chi_q(W_{r+1,a}^{(i)})\chi_q(W_{r-1,aq_i^2}^{(i)})$.\end{lem}

\demo It suffices to check that
$$m_{k,a}^{(i)}m_{k,aq_i^2}^{(i)}A_{i,aq_i^{2r-1}}^{-1}A_{i,aq_i^{2r-3}}^{-1}\cdots A_{i,aq_i}^{-1}$$
$$=k_{\nu(\nu_{r,a}^{(i)})}{\prod}_{\{(j,k,t)|C_{j,i}<0\text{ , }1\leq k\leq -C_{i,j}\text{ and }t\in\{1,\cdots,-C_{j,i}
+E(r_i(r-k)/r_j)\}\}}(k_{\nu(\Lambda_j)}Y_{j,aq_j^{-(2k-1)/C_{i,j}+2(t-1)}}).$$
This follows from a case by case investigation as in the description of lemma \ref{coint}.\qed

\subsubsection{End of the proof of theorem \ref{exact}}

(1) Consequence of corollary \ref{spefus} and lemma \ref{sspe}.

(2) Same proof as for \cite[Theorem 6.1 (2)]{her06}.

(3) To conclude it suffices to prove that $\chi_q(W_{k,a}^{(i)}
*_f
W_{k,aq_i^2}^{(i)})=\chi_q(W_{k,a}^{(i)})\chi_q(W_{k,aq_i^2}^{(i)})$.
This point follows from theorem \ref{lhw}. 

(Observe that we also have a more direct proof of the point (3) from explicit computations : for the type $sl_2$ the result follows from proposition \ref{luhw} and theorem \ref{prod}. For general $\Glie$, the sub $\hat{\U}_i$-module of
$W_{k,a}^{(i)} *_f W_{k,aq_i^2}^{(i)}$ generated by an highest
weight vector is isomorphic to the module obtained in the
$sl_2$-case because $\hat{\U}_i\simeq\U_{q_i}(\hat{sl_2})$ is sub
Hopf algebra of $\U_q(\hat{\Glie})$ for $\Delta_u$. So in
$\chi_q(W_{k,a}^{(i)} *_f W_{k,aq_i^2}^{(i)})$ we have all dominant monomials of
$\chi_q(W_{k,a}^{(i)})\chi_q(W_{k,aq_i^2}^{(i)})$ listed in
\cite[Lemma 5.6]{her06}, and so $E=0$ in the second part of proposition \ref{mainapp}.)\qed

\subsection{Complements, further possible developments and applications}\label{furth} In this section we present questions, possible
developments and applications motivated or related to the results
addressed in this paper. We hope to discuss these points in other
papers.

\subsubsection{Drinfeld coproduct} The algebra $\tilde{\U}_q(\hat{\Glie})$ is defined in section \ref{udelta} (without affine quantum Serre relations (\ref{equadeuxc})). Denote by $\theta_{i,j}^{\pm}(w_1,\cdots,w_s,z)$ the left member of equality (\ref{equadeuxc}) viewed in the algebra $\tilde{\U}_q(\hat{\Glie})$. We conjecture that these elements are quasi-primitive, that it to say that we have in $\tilde{\U}_q'(\hat{\Glie})\hat{\otimes}\tilde{\U}_q'(\hat{\Glie})$ :
$$\Delta_u(\theta_{i,j}^+(w_1,\cdots,w_s,z))=\theta_{i,j}^+(w_1,\cdots,w_s,z)\otimes 1+\phi_i^-(w_1)\cdots\phi_i^-(w_s)\phi_j^-(z)\otimes \theta_{i,j}^+(uw_1,\cdots,uw_s,uz),$$
$$\Delta_u(\theta_{i,j}^-(w_1,\cdots,w_s,z))=1\otimes \theta_{i,j}^-(uw_1,\cdots,uw_s,uz)+\theta_{i,j}^-(w_1,\cdots,w_s,z)\otimes \phi_i^+(uw_1)\cdots\phi_i^+(uw_s)\phi_j^+(uz).$$
The result is known for $C_{i,j}C_{j,i}\leq 3$ (see \cite{e, gr}). This conjecture implies that the Drinfeld coproduct is compatible with affine quantum Serre relations (\ref{equadeuxc}) and that all the results of this paper can be stated for the algebra $\U_q(\hat{\Glie})$ with these relations. In particular the restriction of remark \ref{toref} would be useless.

\subsubsection{$Q$-system and fermionic formulas} Let $\text{Res}$ be the restriction functor from the category of $\U_q(\hat{\Glie})$-modules to the category of $\U_q(\Glie)$-modules. Denote $Q_k^{(i)}=\text{Res}(W_{k,a}^{(i)})$.

We have the following direct consequence of standard results (see \cite{Cha,Cha2}).

\begin{lem} $Q_k^{(i)}$ is well-defined (ie. $\text{Res}(W_{k,a}^{(i)})$ is independent of $a\in\CC^*$).\end{lem}

We give a proof for the convenience of the reader.

\demo Let $\tau_a:\U_q(\hat{\Glie})\rightarrow\U_q(\hat{\Glie})$ be the standard algebra
automorphism defined by $\tau_a(x_{j,m}^{\pm})=a^{\pm
m}x_{j,m}^{\pm}$, $\tau_a(h_{j,r})=a^r h_{j,r}$, $\tau_a(c^{\pm
1/2})=c^{\pm 1/2}$, $\tau_a(k_h)=k_h$ (see \cite{Cha, Cha2}). Then $\U_q(\Glie)$ is
invariant by $\tau_a$ and $W_{k,a}^{(i)}$ is obtained from
$W_{k,1}^{(a)}$ by pull-back by $\tau_a$.\qed

\noindent Consider $\chi$ the usual character morphism for representations in $\text{Mod}(\U_q(\Glie))$ and $\beta$ the linear map on $\Yim$ such that for a monomial $m$, $\beta(m)=\omega(m)$. We have $\beta\circ\chi_q=\chi\circ\text{Res}$ (see \cite[Theorem 3]{Fre}). So theorem \ref{gentsyst} implies :

\begin{thm} Let $i\in I, r\geq 1$. We have :
$$Q_r^{(i)}\otimes Q_r^{(i)}=Q_{r+1}^{(i)}\otimes Q_{r-1}^{(i)}\oplus(\text{Res}(L_{\nu_{r,a}^{(i)}})\otimes({\bigotimes}_{\{(j,k)|C_{j,i}<0\text{ , }1\leq k\leq -C_{i,j}\}}Q_{-C_{j,i}
+E(r_i(r-k)/r_j)}^{(j)})).$$
\end{thm}

Moreover :

\begin{thm}For $i\in I, a\in\CC^*$,
$\mathcal{W}_{k,a}^{(i)}=\chi_q(W_{k,aq_i^{-2k}}^{(i)})Y_{i,aq_i^{-2k}}^{-1}Y_{i,aq_i^{-2k+2}}^{-1}\cdots Y_{i,aq_i^{-2}}^{-1}k_{-k\nu(\Lambda_i^{\vee})}$
considered as a polynomial in $A_{j,b}^{-1}$ has a limit as a formal power series, that is to say 
$$\exists {\text{lim}}_{k\rightarrow \infty}\mathcal{W}_{k,a}^{(i)}\in\ZZ[[A_{j,b}^{-1}]]_{j\in I, b\in\CC^*}.$$

\noindent In particular for $i\in I$, $\mathcal{Q}_k^{(i)}=\chi(Q_k^{(i)})e^{-k\Lambda_i}$ considered as a polynomial in $e^{-\alpha_j}$ has a limit as a formal
power series, that is to say 
$$\exists {\text{lim}}_{k\rightarrow \infty}\mathcal{Q}_k^{(i)}\in\ZZ[[e^{-\alpha_j}]]_{j\in I}.$$\end{thm}

\noindent Such a result was first proved in \cite{Nad} for simply-laced quantum affine algebras and then in \cite{her06} for general quantum affine algebras. These asymptotic properties are a crucial point in the proof of the Kirillov-Reshetikhin conjecture. In particular it seems natural that it should be possible to establish generalized fermionic formulas for the general $Q_k^{(i)}$ as in \cite{hkoty}.

\noindent Besides it should be possible to give an interpretation of
the general $T$-system (theorem \ref{gentsyst}) in terms of
integrable systems as in the original approach \cite{kns}.

\subsubsection{Other possible applications and developments}\label{odev} After an explanation by A. Moura of the proof of \cite{cm}, it seems that the quantum fusion tensor category constructed in this paper will lead to a
description of the block decomposition of the category $\text{Mod}(\U_q(\hat{\Glie}))$ in the spirit of \cite{cm}, the usual tensor product being replaced by the fusion tensor product $*_f$ in the proofs. We hope to treat this point in a separate publication.

\noindent It should be possible to establish a link with the fusion
modules defined in \cite{fl} for classical affine algebras. Indeed we use for $u$ the $\ZZ$-grading of
$\U_q(\hat{\Glie})$ which is the quantum analog of the natural
grading of $\Glie[t^{\pm}]$.

\noindent We also would like to mention the following points. In this paper only the case $q$ generic is considered; an
analog theory for the root of unity cases has still to be
established. As mentioned in the introduction there should be other
examples of the construction used in this paper and a more general
framework involving "quantum Hopf vertex algebras" still to be
defined. Eventually, as we have a forgetful functor to the category
of vector spaces, the tensor category constructed in this paper
should lead to the construction of certain new Hopf algebras $A$ in
a Tannaka-Krein reconstruction spirit.

\end{document}